\documentclass[11pt]{article}

\usepackage{latexsym,amssymb,color}
\usepackage[abbrev]{amsrefs}
\usepackage{amsmath}

\usepackage[truedimen,margin=30truemm, bottom=30truemm]{geometry}

\newtheorem{Theorem}{\bf Theorem}[section]
\newtheorem{Lemma}{\bf Lemma}[section]
\newtheorem{Proposition}{\bf Proposition}[section]
\newtheorem{Corollary}{\bf Corollary}[section]
\newtheorem{Remark}{\bf Remark}[section]
\newtheorem{Example}{\bf Example}[section]
\newtheorem{Definition}{\bf Definition}[section]

\newenvironment{theorem}{\begin{Theorem}$\!\!\!$}{\end{Theorem}}
\newenvironment{lemma}{\begin{Lemma}$\!\!\!$}{\end{Lemma}}

\newenvironment{corollary}{\begin{Corollary}$\!\!\!$}{\end{Corollary}}
\newenvironment{remark}{\begin{Remark}$\!\!\!$}{\end{Remark}}

\newenvironment{definition}{\begin{Definition}$\!\!\!$}{\end{Definition}}

\numberwithin{equation}{section}

\newcommand{\G}{\mathsf{G}}
\newcommand{\dee}{{\rm{d}}}

\begin{document}
\title{Fundamental solution to the heat equation\\
with a dynamical boundary condition}
\author{Kazuhiro Ishige, Sho Katayama, and Tatsuki Kawakami}
\date{}
\maketitle
\begin{abstract}
We give an explicit representation of the fundamental solution to the heat equation on a half-space of ${\mathbb R}^N$ 
with the homogeneous dynamical boundary condition, and obtain upper and lower estimates of the fundamental solution. 
These enable us to obtain sharp decay estimates of solutions to the heat equation with the homogeneous dynamical boundary condition. 
Furthermore, as an application of our decay estimates, 
we identify the so-called Fujita exponent for a semilinear heat equation 
on the half-space of ${\mathbb R}^N$ with the homogeneous dynamical boundary condition. 
\end{abstract}

\vspace{40pt}
\noindent 
Addresses:

\smallskip
\noindent 
K. I.: Graduate School of Mathematical Sciences, The University of Tokyo,\\ 
3-8-1 Komaba, Meguro-ku, Tokyo 153-8914, Japan\\
\noindent 
E-mail: {\tt ishige@ms.u-tokyo.ac.jp}\\

\smallskip
\noindent 
S. K.: Graduate School of Mathematical Sciences, The University of Tokyo,\\ 
3-8-1 Komaba, Meguro-ku, Tokyo 153-8914, Japan\\
\noindent 
E-mail: {\tt katayama-sho572@g.ecc.u-tokyo.ac.jp}\\

\smallskip
\noindent 
{T. K.}: Applied Mathematics and Informatics Course,\\ 
Faculty of Advanced Science and Technology, Ryukoku University,\\
1-5 Yokotani, Seta Oe-cho, Otsu, Shiga 520-2194, Japan\\
\noindent 
E-mail: {\tt kawakami@math.ryukoku.ac.jp}\\

\newpage
%%%%%%%%%%%%%%%%%%%%%%%%%%%%%%%%%%%%%
%%%%%%%%%%%%%%%%%%%%%%%%%%%%%%%%%%%%%
\section{Introduction}
%%%%%%%%%%%%%%%%%%%%%%%%%%%%%%%%%%%%%
%%%%%%%%%%%%%%%%%%%%%%%%%%%%%%%%%%%%%
This paper is concerned with the existence and decay estimates of solutions to 
the heat equation in a half-space of ${\mathbb R}^N$ with the homogeneous dynamical boundary condition
\begin{equation}
\tag{H}
\label{eq:H}
	\left\{
	\begin{array}{ll}
	\partial_t u-\Delta u=0 & \mbox{in}\quad\Omega\times(0,\infty),\vspace{5pt}\\
	\partial_t u+\partial_\nu u=0 & \mbox{on}\quad\partial\Omega\times(0,\infty),\vspace{5pt}\\
	u(x,0)=\phi^i(x) & \mbox{in}\quad\Omega,\vspace{5pt}\\
	u(x,0)=\phi^b(x) & \mbox{on}\quad\partial\Omega,
	\end{array}
	\right.
\end{equation}
where $\phi=(\phi^i,\phi^b)\in L^q(\Omega)\times L^q(\partial\Omega)$ for some $q\in[1,\infty]$ and 
$$
	\Omega:=\{x=(x',x_N)\in{\mathbb R}^N\,:\,x'\in{\mathbb R}^{N-1},\,\,x_N\in(0,\infty)\}. 
$$
Here $N\ge 1$, $\partial_t:=\partial/\partial t$, $\partial_\nu:=\partial/\partial\nu=-\partial/\partial x_N$, 
and $\nu$ is the outer unit normal vector to $\partial\Omega$.  
We often identify $\partial\Omega={\mathbb R}^{N-1}$ if $N\ge 2$ 
and $\Omega=(0,\infty)$ with $\partial\Omega=\{0\}$ if $N=1$. 
In this paper we give an explicit representation of 
the fundamental solution $G$ to problem~\eqref{eq:H} 
and obtain upper and lower pointwise estimates of the fundamental solution~$G$. 
These enable us to establish sharp decay estimates of solutions to problem~\eqref{eq:H}.
Furthermore, as an application of our decay estimates, 
we identify the Fujita exponent for the semilinear heat equation 
with the homogeneous dynamical boundary condition
\begin{equation}
\tag{SH}
\label{eq:SH}
	\left\{
	\begin{array}{ll}
	\displaystyle{\partial_t u=\Delta u+u^p}, & x\in\Omega,\,\,\,t>0,\vspace{5pt}\\
	\displaystyle{\partial_tu+\partial_\nu u=0}, & x\in\partial\Omega,\,\,\, t>0,\vspace{5pt}\\
	\displaystyle{u(x,0)=\phi^i(x)\ge0},\qquad & x\in\Omega,\vspace{5pt}\\
	\displaystyle{u(x,0)=\phi^b(x)\ge0}, & x\in\partial\Omega,
	\end{array}
\right.
\end{equation}
where $p>1$. 
(See Section~3.)

The boundary condition $\partial_t u+\partial_\nu u=0$ on $\partial\Omega\times(0,\infty)$ in problem~\eqref{eq:H} 
describes diffusion through the boundary in processes such as thermal contact with a perfect conductor or
diffusion of solute from a well-stirred fluid or vapour (see e.g., \cite{C}). 
Various aspects of analysis of parabolic equations with dynamical boundary conditions 
have been treated by many authors 
(see e.g., \cites{AF, AQR, BBR, BC, BP01, BP02, DPZ, EMR, E, FIK01, FIK02, FIK02a, FIK02b, 
FIK03, FIK04, FIKL, FQ, FV, HQ, K, PM, R, VV01, VV02, Vi}). 

In \cites{FIK03}
Fila {\it et al.} obtained the existence and local-in-time $L^\infty(\Omega)$-estimates of solutions to problem~\eqref{eq:H} 
with $(\phi^i,\phi^b)\in L^\infty(\Omega)\times L^\infty(\partial\Omega)$ 
by representing solutions via combinations 
of solutions to the Laplace equation with an inhomogeneous dynamical boundary condition 
and an inhomogeneous heat equation with the homogeneous Dirichlet boundary condition 
and by applying the contraction mapping theorem to an integral system involving the Gauss kernel and the Poisson kernel. 
In \cite{FIK04} they developed the arguments in \cite{FIK03} to 
obtain local-in-time decay estimates of solutions to problem~\eqref{eq:H} 
in the case when $\phi^i$ belongs to $L^q$-space in $\Omega$ with some weight, where $q\in[1,\infty]$, and $\phi^b=0$ on $\partial\Omega$.
In \cites{FIK03, FIK04} solutions were represented by a combination of possibly sign-changing functions, 
and it seems difficult to obtain sharp decay estimates of solutions by their arguments. 
In this paper, thanks to our explicit representation of the fundamental solution~$G$, 
we obtain upper and lower pointwise estimates of the fundamental solution~$G$ and 
sharp global-in-time decay estimates of solutions to problem~\eqref{eq:H} 
for $\phi=(\phi^i,\phi^b)\in L^q(\Omega)\times L^q(\partial\Omega)$, where $q\in[1,\infty]$. 
Our arguments improve the results in \cites{FIK03, FIK04} 
and they are quite simpler than those in \cites{FIK03, FIK04}.  
\medskip

We introduce some notation. 
We denote by $\dee\sigma$ the surface measure on $\partial\Omega$ if $N\ge 2$ and $\sigma(\partial\Omega)=\sigma(\{0\})=1$ if $N=1$. 
We set $C_0(\partial\Omega)=C_0({\mathbb R}^{N-1})$ if $N\ge 2$ and $C_0(\partial\Omega)={\mathbb R}$ if $N=1$. 
For any $d=1,2,\dots$, we denote by $\Gamma_d=\Gamma_d(x,t)$ 
the Gauss kernel in ${\mathbb R}^d\times(0,\infty)$, that is,
\begin{equation}
\label{eq:1.1}
	\Gamma_d(x,t):=(4\pi t)^{-\frac{d}{2}}e^{-\frac{|x|^2}{4t}},
	\quad
	(x,t)\in\mathbb R^d\times(0,\infty). 
\end{equation}
Then $\Gamma_d$ satisfies 
\begin{equation}
\label{eq:1.2}
\partial_t \Gamma_d=\Delta\Gamma_d\quad\mbox{in}\quad{\mathbb R}^d\times(0,\infty),
\qquad
\Gamma_d(x,0)=\delta_d(x)\quad\mbox{for}\quad x\in{\mathbb R}^d,
\end{equation}
where $\delta_d$ is the Dirac function in~${\mathbb R}^d$. 

For $N\ge 2$, 
let $p_N$ be the Poisson kernel in $\Omega$, that is, 
$$
p_N(x):=c_Nx_N^{1-N}\left(1+\frac{|x'|^2}{x_N^2}\right)^{-\frac{N}{2}}=c_Nx_N|x|^{-N},
\quad
x=(x',x_N)\in\Omega,
$$
where $c_N$ is a positive constant with 
$$
\int_{{\mathbb R}^{N-1}}p_N(x',x_N)\,\dee x'=1,
\qquad x_N>0. 
$$
For $N=1$, let $p_N(x)=1$ for $x\in\Omega=(0,\infty)$ and $c_N=1$ for convenience. 
Set 
\begin{equation}
\label{eq:1.3}
P_N(x,t):=p_N(x',x_N+t)
=\left\{
\begin{array}{ll}
c_N(x_N+t)|x+te_N|^{-N} & \mbox{if}\quad N\ge 2,\vspace{5pt}\\
1 & \mbox{if}\quad N=1,
\end{array}
\right.\end{equation}
for $x=(x',x_N)\in\Omega$ and $t>0$, where $e_N$ is the unit vector $(0,\dots,0,1)\in{\mathbb R}^N$.
If $N\ge 2$, 
then $P_N$ is the fundamental solution to the Laplace equation in $\Omega$ 
with the homogeneous dynamical boundary condition, that is, $P_N=P_N(x,t)$ satisfies 
$$
\left\{
\begin{array}{ll}
\displaystyle{-\Delta P_N=0} & \mbox{in}\quad \Omega\times(0,\infty),\vspace{5pt}\\
\displaystyle{\partial_t P_N+\partial_\nu P_N=0} & \mbox{on}\quad \partial\Omega\times(0,\infty),\vspace{5pt}\\
\displaystyle{P_N(x',0,0)=\delta_{N-1}(x)} & \mbox{for}\quad x=(x',0)\in\partial\Omega.
\end{array}
\right.
$$

For any $(f,g)\in L^r(\Omega)\times L^r(\partial\Omega)$ with $r\in[1,\infty]$, 
we set
$$
\|(f,g)\|_{L^r(\Omega)\times L^r(\partial\Omega)}:=\|f\|_{L^r(\Omega)}+\|g\|_{L^r(\partial\Omega)}.
$$
Furthermore, for any $f\in C(\overline{\Omega})$, if $f\in L^r(\Omega)$ and $f|_{\partial\Omega}\in L^r(\partial\Omega)$, 
then we write
$$
\|f\|_{L^r(\Omega)\times L^r(\partial\Omega)}:=\|(f,f|_{\partial\Omega})\|_{L^r(\Omega)\times L^r(\partial\Omega)}
$$
for simplicity. Here $f|_{\partial\Omega}$ is the restriction of $f$ on $\partial\Omega$.
\medskip

We state our results on problem~\eqref{eq:H}. 
In Theorem~\ref{Theorem:1.1} 
we give an explicit representation of the fundamental solution to problem~\eqref{eq:H}. 
\begin{theorem}
\label{Theorem:1.1}
Let $N\ge 1$. Set
\begin{equation}
	\label{eq:1.4}
	G(x,y,t):=\Gamma_N(x-y,t)-\Gamma_N(x-y^*,t)+H(x,y,t)
\end{equation}
for $(x,y,t)\in D:=\overline{\Omega}\times\overline{\Omega}\times(0,\infty)$, 
where $y^*=(y',-y_N)$ for $y=(y',y_N)\in\overline{\Omega}$ and 
\begin{equation}
	\label{eq:1.5}
	H(x,y,t):=\,-2\int_0^t (\partial_{x_N}\Gamma_N)(x-y^*+\tau e_N,t-\tau)\,\dee\tau. 
\end{equation}
\begin{itemize}
  \item[{\rm (1)}] 
  $G$ is the fundamental solution to problem~\eqref{eq:H}, that is, 
  \begin{itemize}
  \item[{\rm (a)}]
  for any fixed $y\in\overline\Omega$, 
  the function $G=G(x,y,t)$ satisfies 
  \begin{equation*}
	\left\{
	\begin{array}{ll}
	\displaystyle{\partial_t G-\Delta_x G}=0 
	& \mbox{in}\quad \Omega\times(0,\infty),
	\vspace{5pt}\\
	\displaystyle{\partial_tG+\partial_\nu G=0} 
	& \mbox{on}\quad\partial\Omega\times(0,\infty),
	\end{array}
	\right.
  \end{equation*}
  as a function of the variables $(x,t)\in\overline{\Omega}\times(0,\infty)$;
  \item[{\rm (b)}]
  for any $\zeta\in C_0(\Omega)$, 
  $$
  \lim_{t\to +0}\int_\Omega G(x,y,t)\zeta(x)\,\dee x=\zeta(y),\quad y\in\Omega;
  $$
  \item[{\rm (c)}]
  for any $\eta\in C_0(\partial\Omega)$, 
  $$
  \lim_{t\to +0}\int_{\partial\Omega} G(x,y,t)\eta(x)\,\dee\sigma(x)=\eta(y),\quad y\in\partial\Omega.
  $$
  \end{itemize}
  \item[{\rm (2)}]
  $G(x,y,t)=G(y,x,t)$ for $(x,y,t)\in D$.
  \item[{\rm (3)}]
  For any $(y,t)\in\overline\Omega\times(0,\infty)$, 
  $$
  \int_\Omega G(x,y,t)\,\dee x+\int_{\partial\Omega} G(x,y,t)\,\dee \sigma(x)=1.
  $$
  \end{itemize}
\end{theorem}
In the next theorem we give upper and lower pointwise estimates of the function $H$. 
\begin{theorem}
\label{Theorem:1.2}
There exists $C=C(N)>0$ such that 
\begin{equation}
\label{eq:1.6}
	C^{-1}\underline{H}(x,y,t)\le H(x,y,t)\le C \overline{H}(x,y,t)
\end{equation}
for $(x,y,t)\in D$,
where
\begin{align*}
	 & \overline{H}(x,y,t):=
	\left\{
	\begin{array}{ll}
	\displaystyle{P_N(x-y^*,t)},
	&\quad (x,y,t)\in D_1,
	\vspace{5pt}\\
	\displaystyle{\Gamma_N(x-y^*,t)},
	&\quad (x,y,t)\in D_2\cup D_4,
	\vspace{5pt}\\
	\displaystyle{(x_N+y_N+t)\Gamma_N(x-y^*,t)},\quad
	&\quad (x,y,t)\in D_3,\vspace{2pt}
	\end{array}
	\right.\\
	 & \underline{H}(x,y,t):=
	\left\{
	\begin{array}{ll}
	\displaystyle{P_N(x-y^*,t)},
	&\quad (x,y,t)\in D_1,
	\vspace{5pt}\\
	\displaystyle{\Gamma_N\left(x-y^*,\frac{t}{2}\right)},
	&\quad (x,y,t)\in D_2\cup D_4,
	\vspace{5pt}\\
	\displaystyle{(x_N+y_N+t)\Gamma_N\left(x-y^*,\frac{t}{2}\right)},
	&\quad (x,y,t)\in D_3.
	\end{array}
\right.
\end{align*}
Here
$$
	\begin{aligned}
	&
	D_1:=\{(x,y,t)\in D\,:\,|x-y^*|^2<2(N+2)t,\,\,\,t<4(N+2)\},
	\\
	&
	D_2:=\{(x,y,t)\in D\,:\,|x-y^*|^2<2(N+2)t,\,\,\,t\ge 4(N+2)\},
	\\
	&
	D_3:=\{(x,y,t)\in D\,:\,|x-y^*|^2\ge2(N+2)t,\,\,\, x_N+y_N+t<1\},
	\\
	&
	D_4:=\{(x,y,t)\in D\,:\,|x-y^*|^2\ge2(N+2)t,\,\,\, x_N+y_N+t\ge1\}.
	\end{aligned}
$$
\end{theorem}
Theorem~\ref{Theorem:1.2} reveals the relationship among 
the fundamental solution~$G$ and two kernels $\Gamma_N$ and $P_N$. 

As a corollary of Theorems~\ref{Theorem:1.1} and \ref{Theorem:1.2}, we obtain decay estimates of solutions to problem~\eqref{eq:H} 
for $\phi=(\phi^i,\phi^b)\in L^q(\Omega)\times L^q(\partial\Omega)$ with $q\in[1,\infty]$. 
For any $\phi=(\phi^i,\phi^b)$, where 
$\phi^i$ and $\phi^b$ are measurable functions in $\Omega$ and on $\partial\Omega$, respectively, 
set 
$$	
[\G(t)\phi](x):=\int_\Omega G(x,y,t)\phi^i(y)\,\dee y+\int_{\partial\Omega} G(x,y,t)\phi^b(y)\,\dee \sigma(y)
$$
for $(x,t)\in\overline{\Omega}\times(0,\infty)$ if it is well-defined. 
\begin{corollary}
\label{Corollary:1.1}
\begin{itemize}
\item[{\rm (1)}] 
	Let $\phi\in L^q(\Omega)\times L^q(\partial\Omega)$, where $q\in[1,\infty]$. 
	Then $\G(t)\phi$ is bounded and smooth on $\overline{\Omega}\times[T,\infty)$ for $T>0$
	and it satisfies 
	$$
		\left\{
		\begin{array}{ll}
		\partial_t(\G(t)\phi)-\Delta (\G(t)\phi)=0 & \mbox{in}\quad\Omega\times(0,\infty),\vspace{5pt}\\
		\partial_t(\G(t)\phi)+\partial_\nu(\G(t)\phi)=0 & \mbox{on}\quad\partial\Omega\times(0,\infty),
		\end{array}
		\right.
	$$
	in the classical sense. Furthermore,  
	$$
		\lim_{t\to +0}\|\G(t)\phi-\phi\|_{L^q(\Omega)\times L^q(\partial\Omega)}=0
		\quad\mbox{if}\quad 1\le q<\infty.
	$$
\item[{\rm (2)}] 
  	There exists $C_1=C_1(N)>0$ such that 
  	\begin{equation}
  	\label{eq:1.7}
		\|\G(t)\phi\|_{L^r(\Omega)\times L^r(\partial\Omega)}
		\le C_1\bigg(t^{-\frac{N}{2}\left(\frac{1}{q}-\frac{1}{r}\right)}
		+t^{-(N-1)\left(\frac{1}{q}-\frac{1}{r}\right)}\bigg)\|\phi\|_{L^q(\Omega)\times L^q(\partial\Omega)}
	\end{equation}
  	for $t>0$, $\phi\in L^q(\Omega)\times L^q(\partial\Omega)$, and $1\le q\le r\le\infty$. 
  	In particular, 
  	\begin{equation}
  	\label{eq:1.8}
		\|\G(t)\phi\|_{L^q(\Omega)\times L^q(\partial\Omega)}
		\le\|\phi\|_{L^q(\Omega)\times L^q(\partial\Omega)}
  	\end{equation}
  	for $t>0$, $\phi\in L^q(\Omega)\times L^q(\partial\Omega)$, and $1\le q\le\infty$.
\item[{\rm (3)}]  
	Let $\phi=(\phi^i,\phi^b)$, where $\phi^i$ and $\phi^b$ are nonnegative measurable functions in $\Omega$ and on $\partial\Omega$, respectively. 
  	Then there exists $C_2=C_2(N)>0$ such that 
  	$$
		\begin{aligned}
  		& 
		[\G(t)\phi](x)
		\\
  		& 
		\ge C_2\int_\Omega \Gamma_N\left(x-y^*,\frac{t}{2}\right)\phi^i(y)\,\dee y
		+C_2\int_{\partial\Omega} \Gamma_N\left(x-y^*,\frac{t}{2}\right)\phi^b(y)\,\dee\sigma(y)
  		\end{aligned}
	$$
	for $(x,t)\in \overline{\Omega}\times[1,\infty)$.
\end{itemize}
\end{corollary}
Furthermore, we obtain the semigroup property of $\G(t)$.
\begin{corollary}
\label{Corollary:1.2} 
\begin{itemize}
  \item[{\rm (1)}]
  For any $x$, $y\in\overline{\Omega}$, $t\in[0,\infty)$, and $s\in(0,\infty)$, 
  $$
  G(x,y,t+s)=\int_\Omega G(x,z,t)G(z,y,s)\,\dee z+\int_{\partial\Omega}G(x,z,t)G(z,y,s)\,\dee\sigma(z).
  $$
  \item[{\rm (2)}]
  Let $\phi^i$ and $\phi^b$ be nonnegative measurable functions in $\Omega$ and on $\partial\Omega$, respectively, or 
  $(\phi^i,\phi^b)\in L^q(\Omega)\times L^q(\partial\Omega)$ for some $q\in[1,\infty]$. 
  Set $\phi=(\phi^i,\phi^b)$. Then
  $$
  \G(t)\G(s)\phi=\G(t+s)\phi,\quad t,s\in(0,\infty).
  $$
\end{itemize}
\end{corollary}
As a corollary of Corollaries~\ref{Corollary:1.1} and \ref{Corollary:1.2}, 
we show that decay estimates in Corollary~\ref{Corollary:1.1}-(2) are sharp. 
\begin{corollary}
\label{Corollary:1.3} 
Let $1\le q\le r\le\infty$. Set 
$$
\|\G(t)\|_{q\to r}:=\sup\left\{\frac{\|\G(t)\phi\|_{L^r(\Omega)\times L^r(\partial\Omega)}}{\|\phi\|_{L^q(\Omega)\times L^q(\partial\Omega)}}\,:\,
\phi=(\phi^i,\phi^b)\in L^q(\Omega)\times L^q(\partial\Omega),\,\,\phi\not=(0,0)\right\}.
$$
Then there exists $C=C(N)>0$ such that 
\begin{align*}
C^{-1}t^{-(N-1)\left(\frac{1}{q}-\frac{1}{r}\right)}
 & \le \|\G(t)\|_{q\to r}\le Ct^{-(N-1)\left(\frac{1}{q}-\frac{1}{r}\right)},\quad t\in(0,1],\\
C^{-1}t^{-\frac{N}{2}\left(\frac{1}{q}-\frac{1}{r}\right)}
 & \le \|\G(t)\|_{q\to r}\le Ct^{-\frac{N}{2}\left(\frac{1}{q}-\frac{1}{r}\right)},\qquad\,\, t\in(1,\infty).
\end{align*}
\end{corollary}
\vspace{5pt}

The rest of this paper is organized as follows. 
In Section~2 we prove Theorems~\ref{Theorem:1.1} and \ref{Theorem:1.2}. 
Furthermore, we prove Corollaries~\ref{Corollary:1.1}, \ref{Corollary:1.2}, and \ref{Corollary:1.3}. 
In Section~3 we apply our results on problem~\eqref{eq:H}  to identify the so-called Fujita exponent 
for problem~\eqref{eq:SH}. 
%%%%%%%%%%%%%%%%%%%%%%%%%%%%%%%%%%%%%
%%%%%%%%%%%%%%%%%%%%%%%%%%%%%%%%%%%%%
\section{Proof of Theorems~\ref{Theorem:1.1} and \ref{Theorem:1.2}}
%%%%%%%%%%%%%%%%%%%%%%%%%%%%%%%%%%%%%
%%%%%%%%%%%%%%%%%%%%%%%%%%%%%%%%%%%%%
In all that follows we will use $C$ to denote generic positive constants which are independent of $(x,y,t)\in D$ 
and point out that $C$  may take different values  within a calculation.
\vspace{5pt}
\newline
{\bf Proof of Theorem~\ref{Theorem:1.1}.}
We observe from \eqref{eq:1.2} and the definition of $H$ (see \eqref{eq:1.5}) that, 
for any fixed $y\in\overline{\Omega}$, $G$ satisfies 
$\partial_t G-\Delta_x G=0$
as a function of the variable $(x,t)\in\Omega\times(0,\infty)$. 
Furthermore, since 
\begin{equation}
\label{eq:2.1}
	\Gamma_N(x-y,t)-\Gamma_N(x-y^*,t)
	=(4\pi t)^{-\frac{N}{2}}e^{-\frac{|x'-y'|^2}{4t}}
	\left(e^{-\frac{|x_N-y_N|^2}{4t}}-e^{-\frac{|x_N+y_N|^2}{4t}}\right),
\end{equation}
we see that
\begin{align*}
	& 
	\partial_t H(x,y,t)+\partial_\nu H(x,y,t)
	\\
 	& 
	=-2\int_0^t (\partial_t-\partial_{x_N})(\partial_{x_N}\Gamma_N)(x-y^*+\tau e_N,t-\tau)\,\dee\tau
	\\
 	& 
	=2\int_0^t \frac{\partial}{\partial\tau}(\partial_{x_N}\Gamma_N)(x-y^*+\tau e_N,t-\tau)\,\dee\tau
	\\
 	& 
	=-2(\partial_{x_N}\Gamma_N)(x'-y',y_N,t)
 	=-(\partial_t+\partial_\nu)\left(\Gamma_N(x-y,t)-\Gamma_N(x-y^*,t)\right)
\end{align*}
for $x\in\partial\Omega$, $y\in\overline{\Omega}$, and $t>0$.
This implies that $G$ satisfies 
$\partial_t G+\partial_\nu G=0$
as a function of the variables $(x,t)\in\partial\Omega\times(0,\infty)$. 
In addition, it follows from \eqref{eq:1.5} that  
\begin{equation}
	\label{eq:2.2}
	\begin{aligned}
	H(x,y,t)= 
	& 
	-2\int_0^t (\partial_{x_N}\Gamma_N)(x-y^*+\tau e_N,t-\tau)\,\dee\tau
	\\
	= & \int_0^t(4\pi(t-\tau))^{-\frac{N}{2}}\frac{x_N+y_N+\tau}{t-\tau}
	e^{-\frac{|x-y^*+\tau e_N|^2}{4(t-\tau)}}\,\dee\tau
	\\
	= & 
	\int_0^t(4\pi(t-\tau))^{-\frac{1}{2}}\frac{x_N+y_N+\tau}{t-\tau}
	e^{-\frac{|x_N+y_N+\tau|^2}{4(t-\tau)}}\Gamma_{N-1}(x'-y',t-\tau)\,\dee\tau
	\end{aligned}
\end{equation}
for $(x,y,t)\in D$.
This together with \eqref{eq:2.1} implies assertion~(2).

We prove assertion~(3). 
Let $(y,t)\in\overline{\Omega}\times(0,\infty)$. 
Since $\partial\Omega={\mathbb R}^{N-1}\times\{0\}$ and 
$$
	\int_{{\mathbb R}^{N-1}}\Gamma_{N-1}(x',t)\,\dee x'=1,\quad t>0,
$$ 
by \eqref{eq:2.2} we have 
\begin{align*}
	\int_{\partial\Omega}H(x,y,t)\,\dee\sigma(x)
	&
	=\int_0^t(4\pi(t-\tau))^{-\frac{1}{2}}\frac{y_N+\tau}{t-\tau}
	e^{-\frac{(y_N+\tau)^2}{4(t-\tau)}}\,\dee\tau
	\\
	&
	=\int_0^t(\pi \tau)^{-\frac{1}{2}}\frac{y_N+t-\tau}{2\tau}
	e^{-\frac{(y_N+t-\tau)^2}{4\tau}}\,\dee\tau,\\
	\int_\Omega H(x,y,t)\,\dee x
	&
	=
	\int_0^t\int_0^\infty(4\pi(t-\tau))^{-\frac{1}{2}}\frac{x_N+y_N+\tau}{t-\tau}
	e^{-\frac{(x_N+y_N+\tau)^2}{4(t-\tau)}}\,\dee x_N\,\dee\tau
	\\
	&
	=\int_0^t(\pi(t-\tau))^{-\frac{1}{2}}
	\int_{\frac{y_N+\tau}{2(t-\tau)^{1/2}}}^\infty 2\eta e^{-\eta^2}\,\dee\eta\,\dee\tau
	\\
	&
	=\int_0^t(\pi (t-\tau))^{-\frac{1}{2}}
	e^{-\frac{(y_N+\tau)^2}{4(t-\tau)}}\,\dee\tau
	=\int_0^t(\pi \tau)^{-\frac{1}{2}}
	e^{-\frac{(y_N+t-\tau)^2}{4\tau}}\,\dee\tau.
\end{align*}
These imply that 
\begin{equation}
\label{eq:2.3}
	\begin{aligned}
	&
	\int_\Omega H(x,y,t)\,\dee x+\int_{\partial\Omega}H(x,y,t)\,\dee\sigma(x)
	\\
	&
	=
	\int_0^t(\pi \tau)^{-\frac{1}{2}}\bigg(1+\frac{y_N+t-\tau}{2\tau}\bigg)
	e^{-\frac{(y_N+t-\tau)^2}{4\tau}}\,\dee\tau
	=2\pi^{-\frac{1}{2}}\int_{\frac{y_N}{2t^{1/2}}}^\infty e^{-\xi^2}\, \dee\xi
	\end{aligned}
\end{equation}
and 
\begin{equation}
\label{eq:2.4}
	\lim_{t\to+0}\int_{\partial\Omega}H(x,y,t)\,\dee\sigma(x)=0\quad\mbox{for $y\in\Omega$},\quad
	\lim_{t\to+0}\int_\Omega H(x,y,t)\,\dee x=0\quad\mbox{for $y\in\partial\Omega$}.
\end{equation}
On the other hand, it follows from \eqref{eq:2.1} that
\begin{equation}
\label{eq:2.5}
	\begin{aligned}
	\int_\Omega \bigg(\Gamma_N(x-y,t)-\Gamma_N(x-y^*,t)\bigg)\,\dee x
	 & =\int_0^\infty (4\pi t)^{-\frac{1}{2}}
	\bigg(e^{-\frac{(x_N-y_N)^2}{4t}}-e^{-\frac{(x_N+y_N)^2}{4t}}\bigg)\,\dee x_N
	\\
	&
	=\pi^{-\frac{1}{2}}\int_{-\frac{y_N}{2t^{1/2}}}^{\frac{y_N}{2t^{1/2}}}e^{-\xi^2}\,\dee\xi
	=2\pi^{-\frac{1}{2}}\int_0^{\frac{y_N}{2t^{1/2}}}e^{-\xi^2}\,\dee\xi.
	\end{aligned}
\end{equation}
Since $\Gamma_N(x-y,t)=\Gamma_N(x-y^*,t)$ for $x\in\partial\Omega$,
it follows from \eqref{eq:1.4}, \eqref{eq:2.3}, and \eqref{eq:2.5} that
\begin{equation}
\label{eq:2.6}
	\begin{aligned}
	&
	\int_\Omega G(x,y,t)\,\dee x+\int_{\partial\Omega}G(x,y,t)\,\dee\sigma(x)
	\\
	&
	=\int_\Omega \bigg(\Gamma_N(x-y,t)-\Gamma_N(x-y^*,t)\bigg)\,\dee x
	+\int_\Omega H(x,y,t)\,\dee x+\int_{\partial\Omega}H(x,y,t)\,\dee\sigma(x)
	\\
	&
	=2\pi^{-\frac{1}{2}}\int_0^{\frac{y_N}{2t^{1/2}}}e^{-\xi^2}\,\dee\xi
	+2\pi^{-\frac{1}{2}}\int_{\frac{y_N}{2t^{1/2}}}^\infty e^{-\xi^2}\, d\xi
	=2\pi^{-\frac{1}{2}}\int_0^\infty e^{-\xi^2}\,\dee\xi
	=1
	\end{aligned}
\end{equation}
for $(y,t)\in\overline{\Omega}\times(0,\infty)$. 
Thus assertion~(3) holds. 
Similarly, we see that 
\begin{equation}
\label{eq:2.7}
	\lim_{t\to +0}\int_{\Omega\setminus B(y,\delta)} G(x,y,t)\,\dee x=0,
	\quad
	\lim_{t\to +0}\int_{\partial\Omega\setminus B(y,\delta)}G(x,y,t)\,\dee\sigma(x)=0,
\end{equation}
for any $y\in\overline{\Omega}$ and $\delta>0$.
Furthermore, it follows from \eqref{eq:2.4} that
\begin{equation}
\label{eq:2.8}
	\lim_{t\to+0}\int_{\partial\Omega}G(x,y,t)\,\dee\sigma(x)=0\quad\mbox{for $y\in\Omega$},\quad
	\lim_{t\to+0}\int_\Omega G(x,y,t)\,\dee x=0\quad\mbox{for $y\in\partial\Omega$}.
\end{equation}

It remains to prove assertions~(1)-(b) and (1)-(c). 
Let $y\in\Omega$ and $\zeta\in C_0(\Omega)$. 
Then it follows from \eqref{eq:2.6} that
\begin{align*}
 	& 
	\left|\int_\Omega G(x,y,t)\zeta(x)\,\dee x-\zeta(y)\right|
	\\
 	& 
	=\left|-\zeta(y)\int_{\partial\Omega} G(x,y,t)\,\dee \sigma(x)+\int_\Omega G(x,y,t)(\zeta(x)-\zeta(y))\,\dee x\right|
	\\
 	& 
	\le|\zeta(y)|\int_{\partial\Omega} G(x,y,t)\,\dee \sigma(x)
 	+2\|\zeta\|_{L^\infty(\Omega)}\int_{\Omega\setminus B(y,\delta)} G(x,y,t)\,\dee x
	\\
 	& 
	\qquad\quad
 	+\sup_{x\in B(y,\delta)}|\zeta(x)-\zeta(y)|\int_{\Omega\cap B(y,\delta)} G(x,y,t)\,\dee x
\end{align*}
for $\delta>0$.
Then, by \eqref{eq:2.6}, \eqref{eq:2.7}, and \eqref{eq:2.8} we have
$$
	\limsup_{t\to +0}\left|\int_\Omega G(x,y,t)\zeta(x)\,\dee x-\zeta(y)\right|
	\le \sup_{x\in B(y,\delta)}|\zeta(x)-\zeta(y)|.
$$
Since $\delta$ is arbitrary, the continuity of $\zeta$ implies that
$$
	\lim_{t\to +0}\int_\Omega G(x,y,t)\zeta(x)\,\dee x=\zeta(y).
$$
Similarly, if $N\ge 2$, then, 
for any $y\in\partial\Omega$ and $\eta\in C_0(\partial\Omega)$, 
since
\begin{align*}
 	& 
	\left|\int_{\partial\Omega} G(x,y,t)\eta(x)\,\dee \sigma(x)-\eta(y)\right|
	\\
 	& 
	=\left|-\eta(y)\int_\Omega G(x,y,t)\,\dee x+\int_{\partial\Omega} G(x,y,t)(\zeta(x)-\zeta(y))\,\dee\sigma(x)\right|
	\\
 	& 
	\le|\eta(y)|\int_\Omega G(x,y,t)\,\dee x
 	+2\|\eta\|_{L^\infty(\Omega)}\int_{\partial\Omega\setminus B(y,\delta)} G(x,y,t)\,\dee\sigma(x)
	\\
 	& 
	\qquad\quad
 	+\sup_{x\in B(y,\delta)}|\eta(y)-\eta(x)|\int_{\partial\Omega\cap B(y,\delta)} G(x,y,t)\,\dee\sigma(x),
\end{align*}
we have 
$$
	\limsup_{t\to +0}\left|\int_{\partial\Omega} G(x,y,t)\eta(x)\,\dee\sigma(x)-\eta(y)\right|
	\le \sup_{x\in B(y,\delta)}|\eta(x)-\eta(y)|.
$$
Since $\delta$ is arbitrary, the continuity of $\eta$ implies that
$$
	\lim_{t\to +0}\int_{\partial\Omega} G(x,y,t)\eta(x)\,\dee x=\eta(y)\quad\mbox{if}\quad N\ge 2. 
$$
If $N=1$, then, for any $\eta\in C_0(\partial\Omega)={\mathbb R}$ and $y\in\partial\Omega=\{0\}$,
it follows from \eqref{eq:2.6} and \eqref{eq:2.8} that
$$
	\lim_{t\to +0}\int_{\partial\Omega} G(x,y,t)\eta\,\dee\sigma(x)=\eta.
$$
Thus assertions~(1)-(b) and (1)-(c) hold. The proof of Theorem~\ref{Theorem:1.1} is complete.
$\Box$ \vspace{5pt}

\noindent
{\bf Proof of Theorem~\ref{Theorem:1.2}.}
The proof is divided into several steps. 
\vspace{5pt}
\newline
\underline{Step 1.} 
Let $(x,y,t)\in D_1$, and set
$z:=x-y^*+te_N$. 
Since 
$$
	|x-y^*+\tau e_N|^2=|x-y^*+te_N-(t-\tau)e_N|^2=|z|^2-2(t-\tau)z_N+(t-\tau)^2,
$$
by \eqref{eq:2.2} we have 
\begin{equation}
\label{eq:2.9}
	\begin{aligned}
	H(x,y,t)
	 & =\int_0^t(4\pi(t-\tau))^{-\frac{N}{2}}\frac{x_N+y_N+\tau}{t-\tau}
	e^{-\frac{|x-y^*+\tau e_N|^2}{4(t-\tau)}}\,\dee\tau
	\\
	&
	=\int_0^t(4\pi(t-\tau))^{-\frac{N}{2}}\frac{z_N-t+\tau}{t-\tau}
	e^{-\frac{|z|^2-2(t-\tau)z_N+(t-\tau)^2}{4(t-\tau)}}\,\dee\tau
	\\
	&
	=e^{\frac{z_N}{2}}\int_0^t(4\pi(t-\tau))^{-\frac{N}{2}}\bigg(\frac{z_N}{t-\tau}-1\bigg)
	e^{-\frac{|z|^2}{4(t-\tau)}-\frac{t-\tau}{4}}\,\dee\tau
	\\
	&
	=(4\pi)^{-\frac{N}{2}}e^{\frac{z_N}{2}}\int_0^{\frac{t}{|z|^2}}(|z|^2\eta)^{-\frac{N}{2}}
	\bigg(\frac{z_N}{|z|^2\eta}-1\bigg)e^{-\frac{1}{4\eta}-\frac{|z|^2\eta}{4}}|z|^2\,\dee\eta
	\\
	&
	=(4\pi)^{-\frac{N}{2}}z_N|z|^{-N}e^{\frac{z_N}{2}}
	\int_0^{\frac{t}{|z|^2}}
	\bigg(1-z_N^{-1}|z|^2\eta\bigg)\eta^{-\frac{N+2}{2}}e^{-\frac{1}{4\eta}-\frac{|z|^2\eta}{4}}\,\dee\eta.
	\end{aligned}
\end{equation}
Since $z_N\le |x-y^*|+t\le \sqrt{2(N+2)t}+t\le C$, 
by \eqref{eq:1.3} and \eqref{eq:2.9} we obtain
\begin{equation}
\label{eq:2.10}
	H(x,y,t)
	\le C z_N|z|^{-N}e^{\frac{z_N}{2}}\int_0^\infty
	\eta^{-\frac{N+2}{2}}e^{-\frac{1}{4\eta}}\,\dee\eta
	\le Cz_N|z|^{-N}
	=C\overline{H}(x,y,t).
\end{equation}
On the other hand, it follows that
$$
	z_N\ge t,\qquad
	1-z_N^{-1}|z|^2\eta\ge\frac{1}{2}\quad\mbox{for}\quad\eta\in\bigg(0,\frac{t}{2|z|^2}\bigg).
$$
These together with \eqref{eq:2.9} imply that
\begin{equation}
\label{eq:2.11}
\begin{split}
	H(x,y,t)
	&
	\ge
	Cz_N|z|^{-N}e^{\frac{z_N}{2}}
	\int_0^{\frac{t}{2|z|^2}}
	\eta^{-\frac{N+2}{2}}e^{-\frac{1}{4\eta}-\frac{|z|^2\eta}{4}}\,\dee\eta
	\\
	&
	\ge
	Cz_N|z|^{-N}e^{\frac{z_N}{2}-\frac{t}{8}}
	\int_0^{\frac{t}{2|z|^2}}
	\eta^{-\frac{N+2}{2}}e^{-\frac{1}{4\eta}}\,\dee\eta
	\\
	&
	\ge
	Cz_N|z|^{-N}
	\int_0^{\frac{t}{2|z|^2}}
	\eta^{-\frac{N+2}{2}}e^{-\frac{1}{4\eta}}\,\dee\eta.
\end{split}
\end{equation}
Since
$$
	\frac{|z|^2}{t}\le 2\frac{|x-y^*|^2+t^2}{t}
	\le 4(N+2)+2t\le C,
$$
it follows from \eqref{eq:2.11} that $H(x,y,t)\ge C\underline{H}(x,y,t)$ for $(x,y,t)\in D_1$. 
This together with \eqref{eq:2.10} implies \eqref{eq:1.6} for $(x,y,t)\in D_1$. 
\vspace{5pt}
\newline
\underline{Step 2.} 
Let $(x,y,t)\in D_2$.
We find a unique $\tau_*\in(0,t/2]$ such that
\begin{equation}
\label{eq:2.12}
	t-\tau_*=\frac{|x-y^*+\tau_*e_N|^2}{2(N+2)}.
\end{equation}
Indeed, setting
$$
f(\tau):=\frac{|x-y^*+\tau e_N|^2}{2(N+2)}-(t-\tau),
$$
we see that $f(0)<0$ by $(x,y,t)\in D_2$ and $f$ is strictly increasing in $[0,\infty)$.
Furthermore, since $x_N+y_N\ge0$, we have
$$
	f\left(\frac{t}{2}\right)\ge \frac{(x_N+y_N+t/2)^2}{2(N+2)}-\frac{t}{2}
	\ge \frac{t^2}{8(N+2)}-\frac{t}{2}\ge0
$$
for $t\ge 4(N+2)$.
These together with the intermediate value theorem imply the unique existence of 
$\tau_*\in(0,t/2]$ satisfying \eqref{eq:2.12}. 

It follows from \eqref{eq:2.12} that 
\begin{equation}
\label{eq:2.13}
	\begin{aligned}
	&
	\frac{|x-y^*+\tau e_N|^2}{2(N+2)}\le t-\tau_*\le t-\tau,\qquad \tau\in(0,\tau_*],
	\\
	&
	t-\tau\le t-\tau_*\le\frac{|x-y^*+\tau e_N|^2}{2(N+2)},\qquad \tau\in[\tau_*,t).
	\end{aligned}
\end{equation}
On the one hand, since 
$$
	(\partial_t\partial_{x_N}\Gamma_N)(x,t)
	=\left(N+2-\frac{|x|^2}{2t}\right)\frac{x_N}{4t^2}\Gamma_N(x,t),
	\quad (x,t)\in\Omega\times(0,\infty), 
$$
we see that 
\begin{equation}
\label{eq:2.14}
	\mbox{$-(\partial_{x_N}\Gamma_N)(x,\cdot)$ is }
	\left\{
	\begin{array}{l}
	\mbox{increasing in $\displaystyle{\left(0, \frac{|x|^2}{2(N+2)}\right)}$,}\vspace{8pt}\\
	\mbox{decreasing in $\displaystyle{\left(\frac{|x|^2}{2(N+2)},\infty\right)}$.}
	\end{array}
	\right.
\end{equation}
By \eqref{eq:2.13} and \eqref{eq:2.14} we see that 
$$
-\partial_{x_N}\Gamma(x-y_*+\tau e_N,t-\tau)\le -\partial_{x_N}\Gamma(x-y_*+\tau e_N,t-\tau_*)
$$
for $\tau\in(0,t)$, and hence we obtain
\begin{equation}
\label{eq:2.15}
	\begin{aligned}
	H(x,y,t)
	&
	\le -2\int_0^t(\partial_{x_N}\Gamma_N)(x-y^*+\tau e_N,t-\tau_*)\, \dee\tau
	\\
	& 
	=-2\int_0^t \frac{\partial}{\partial\tau}\Gamma_N(x-y^*+\tau e_N,t-\tau_*)\, \dee\tau
	\\
	&
	\le 2\Gamma_N(x-y^*,t-\tau_*)
	\le 2\left(\frac{t}{t-\tau_*}\right)^{\frac{N}{2}}\Gamma_N(x-y^*,t)
	\le C\Gamma_N(x-y^*,t).
	\end{aligned}
\end{equation}

On the other hand,
since 
$$
	\begin{aligned}
	-(\partial_{x_N}\Gamma_N)(x,t-\tau)
	 & 
	 =(4\pi(t-\tau))^{-\frac{N}{2}}\frac{x_N}{2(t-\tau)}e^{-\frac{|x|^2}{4(t-\tau)}}
	\\
	& 
	= \left(\frac{t}{t-\tau}\right)^{\frac{N+2}{2}}
	(4\pi t)^{-\frac{N}{2}}\frac{x_N}{2t}e^{-\frac{|x|^2}{4(t-\tau)}}
	\\
	& 
	\ge (4\pi t)^{-\frac{N}{2}}\frac{x_N}{2t}e^{-\frac{|x|^2}{2t}}
	= -2^{-\frac{N+2}{2}}(\partial_{x_N}\Gamma_N)\left(x,\frac{t}{2}\right)
	\end{aligned}
$$
for $x\in \mathbb R^N$ and $0<\tau\le t/2$,
we have
\begin{equation}
\label{eq:2.16}
	\begin{aligned}
	H(x,y,t)
	&
	\ge
	-2^{-\frac{N+2}{2}}\int_0^{t/2}(\partial_{x_N}\Gamma_N)
	\left(x-y^*+\tau e_N,\frac{t}{2}\right)\, \dee\tau
	\\
	& 
	=-2^{-\frac{N+2}{2}}\int_0^{t/2}\frac{\partial}{\partial\tau}
	\Gamma_N\left(x-y^*+\tau e_N,\frac{t}{2}\right)\, \dee\tau
	\\
	&
	=2^{-\frac{N+2}{2}}\bigg(\Gamma_N\left(x-y^*,\frac{t}{2}\right)-\Gamma_N
	\left(x-y^*+\frac{t}{2}e_N,\frac{t}{2}\right)\bigg)
	\\
	&
	=2^{-\frac{N+2}{2}}\left(1-e^{-\frac{4(x_N+y_N)+t}{8}}\right)\Gamma_N\left(x-y^*,\frac{t}{2}\right).
	\end{aligned}
\end{equation}
Since $x_N+y_N\ge0$ and $t\ge 4(N+2)$ by $(x,y,t)\in D_2$, 
we observe from \eqref{eq:2.16}  that 
$$
	H(x,y,t)\ge C\Gamma_N\left(x-y^*,\frac{t}{2}\right).
$$
This together with \eqref{eq:2.15} implies \eqref{eq:1.6} for $(x,y,t)\in D_2$.
\vspace{5pt}
\newline
\underline{Step 3.} 
Let $(x,y,t)\in D_3\cup D_4$. 
Since 
$$
	t-\tau\le t\le \frac{|x-y^*|^2}{2(N+2)}\le  \frac{|x-y^*+\tau|^2}{2(N+2)},\qquad 0\le \tau<t,
$$
by \eqref{eq:2.14} we have 
$$
	-(\partial_{x_N}\Gamma_N)(x-y^*+\tau e_N,t-\tau)
	\le-(\partial_{x_N}\Gamma_N)(x-y^*+\tau e_N,t).
$$
Then, applying the same arguments as in \eqref{eq:2.15} and \eqref{eq:2.16}, 
we see that
$$
	\begin{aligned}
	H(x,y,t)
	&
	\le -2\int_0^t(\partial_{x_N}\Gamma_N)(x-y^*+\tau e_N,t)\, \dee\tau
	\\
	& 
	=-2\int_0^t \frac{\partial}{\partial\tau}\Gamma_N(x-y^*+\tau e_N,t)\, \dee\tau
	\\
	& 
	\le 2(\Gamma_N(x-y^*,t)-\Gamma_N(x-y^*+te_N,t))
	\\
	& 
	=2\Gamma_N(x-y^*)\left(1-e^{-\frac{2(x_N+y_N)+t}{4}}\right)
	\\
	& 
	\le
	\left\{
	\begin{array}{ll}
	C\Gamma_N(x-y^*,t)(x_N+y_N+t) & \mbox{if}\quad x_N+y_N+t<1,\vspace{5pt}\\
	2\Gamma_N(x-y^*,t) & \mbox{if}\quad x_N+y_N+t\ge 1,
	\end{array}
	\right.
	\end{aligned}
$$
and
$$ 
	\begin{aligned}
	H(x,y,t) 
	& 
	\ge 2^{-\frac{N+2}{2}}\left(1-e^{-\frac{4(x_N+y_N)+t}{8}}\right)\Gamma_N\left(x-y^*,\frac{t}{2}\right)
	\\
 	& 
	\ge\left\{
	\begin{array}{ll}
	\displaystyle{C(x_N+y_N+t)}\Gamma_N\left(x-y^*,\frac{t}{2}\right) 
	& \mbox{if}\quad x_N+y_N+t<1,\vspace{5pt}\\
	\displaystyle{C\Gamma_N\left(x-y^*,\frac{t}{2}\right)} & \mbox{if}\quad x_N+y_N+t\ge 1.
	\end{array}
	\right.
	\end{aligned}
$$
Thus \eqref{eq:1.6} holds for $(x,y,t)\in D$, 
and the proof of Theorem~\ref{Theorem:1.2} is complete. 
$\Box$
\vspace{5pt}

\noindent
{\bf Proof of Corollary~\ref{Corollary:1.1}.}
Let $1\le q\le\infty$ and $\phi=(\phi^i,\phi^b)\in L^q(\Omega)\times L^q(\partial\Omega)$. 
If $1\le q<\infty$, by H\"older's inequality and Theorem~\ref{Theorem:1.1}-(3) we have
\begin{equation}
\label{eq:2.17}
\begin{split}
	 & \sup_{x\in\overline{\Omega}}|[\G(t)\phi](x)| \\
	& 
	\le \sup_{x\in\overline{\Omega}}\left(\int_\Omega G(x,y,t)|\phi^i(y)|^q\,\dee y\right)^{\frac{1}{q}}
	+\sup_{x\in\overline{\Omega}}\left(\int_{\partial\Omega} G(x,y,t)|\phi^b(y)|^q\,\dee\sigma(y)\right)^{\frac{1}{q}}
	\\
 	& 
	\le\sup_{x\in\overline{\Omega}}\|G(x,\cdot,t)\|_{L^\infty(\Omega)}^{\frac{1}{q}}
	(\|\phi^i\|_{L^q(\Omega)}+\|\phi^b\|_{L^q(\partial\Omega)})<\infty
\end{split}
\end{equation}
for $t>0$. This inequality also holds in the case of $q=\infty$. 
Then, by Theorem~\ref{Theorem:1.1} we apply the standard theory for convolutions to obtain assertion~(1).

We prove assertion~(2). 
Since assertion~(2) with $q=\infty$ easily follows from Theorem~\ref{Theorem:1.1}-(3), 
we treat the case of $1\le q<\infty$. 
If $q=1$, by Theorem~\ref{Theorem:1.1}-(3) we obtain
\begin{align*}
 & \|\G(t)\phi\|_{L^1(\Omega)}+\|\G(t)\phi\|_{L^1(\partial\Omega)}\\
 & \le \int_\Omega\int_\Omega G(x,y,t)|\phi^i(y)|\,\dee y\,\dee x+\int_\Omega\int_{\partial\Omega}G(x,y,t)|\phi^b(y)|\,\dee\sigma(y)\,\dee x\\
 & \qquad
 + \int_{\partial\Omega}\int_\Omega G(x,y,t)|\phi^i(y)|\,\dee y\,\dee\sigma(x)+\int_{\partial\Omega}\int_{\partial\Omega}G(x,y,t)|\phi^b(y)|\,\dee\sigma(y)\,\dee\sigma(x)\\
 & \le\int_\Omega\left(\int_\Omega G(x,y,t)\,\dee x+\int_{\partial\Omega}G(x,y,t)\,\dee\sigma(x)\right)|\phi^i(y)|\,\dee y\\
 & \qquad
 +\int_{\partial\Omega}\left(\int_\Omega G(x,y,t)\,\dee x+\int_{\partial\Omega}G(x,y,t)\,\dee\sigma(x)\right)|\phi^b(y)|\,\dee \sigma(y)\\
 & =\|\phi^i\|_{L^1(\Omega)}+\|\phi^b\|_{L^1(\partial\Omega)}=\|\phi\|_{L^1(\Omega)\times L^1(\partial\Omega)}.
\end{align*}
This implies that \eqref{eq:1.8} holds for $q=1$. 
Since \eqref{eq:1.8} holds for $q=\infty$, 
we apply the Riesz-Thorin interpolation theorem to obtain \eqref{eq:1.8} for $q\in[1,\infty]$.

On the other hand, 
since $P_N(x-y^*,t)\le c_N(x_N+y_N+t)^{-(N-1)}\le c_Nt^{-(N-1)}$ for $(x,y,t)\in D$ (see \eqref{eq:1.3}), 
it follows from Theorem~\ref{Theorem:1.2} and \eqref{eq:1.4} that 
$$
G(x,y,t)\le Ct^{-\frac{N}{2}}+C\Gamma_N(x-y^*,t)+CP_N(x-y^*,t)\le Ct^{-\frac{N}{2}}+Ct^{-(N-1)},\quad (x,y,t)\in D. 
$$
This together with \eqref{eq:2.17} implies that 
\begin{equation}
\label{eq:2.18}
	\|\G(t)\phi\|_{L^\infty(\Omega)}+\|\G(t)\phi\|_{L^\infty(\partial\Omega)}
	\le C\left(t^{-\frac{N}{2q}}+t^{-\frac{N-1}{q}}\right)\|\phi\|_{L^q(\Omega)\times L^q(\partial\Omega)}
\end{equation}
for $t>0$. By \eqref{eq:1.8} and \eqref{eq:2.18}, 
for any $r\in[q,\infty]$, we obtain
$$
	\begin{aligned}
	\|\G(t)\phi\|_{L^r(\Omega)\times L^r(\partial\Omega)}
	 & =\|\G(t)\phi\|_{L^r(\Omega)}+\|\G(t)\phi\|_{L^r(\partial\Omega)}
	\\
	& 
	\le \|\G(t)\phi\|_{L^\infty(\Omega)}^{1-\frac{q}{r}}\|\G(t)\phi\|_{L^q(\Omega)}^{\frac{q}{r}}
	+\|\G(t)\phi\|_{L^\infty(\partial\Omega)}^{1-\frac{q}{r}}\|\G(t)\phi\|_{L^q(\partial\Omega)}^{\frac{q}{r}}
	\\
	 & 
	 \le C\left(t^{-\frac{N}{2}\left(\frac{1}{q}-\frac{1}{r}\right)}
	 +t^{-(N-1)\left(\frac{1}{q}-\frac{1}{r}\right)}\right)\|\phi\|_{L^q(\Omega)\times L^q(\partial\Omega)}
	\end{aligned}
$$
for $t>0$. This implies \eqref{eq:1.7}. Thus assertion~(2) follows.

We prove assertion~(3). 
Let $\phi^i$ and $\phi^b$ be nonnegative measurable functions in $\Omega$ and $\partial\Omega$, respectively. 
Since $\Gamma_N(x-y,t)-\Gamma_N(x-y^*,t)$ is nonnegative for $(x,y,t)\in D$ and 
\begin{equation}
\label{eq:2.19}
\{y\in\overline{\Omega}\,:\,(x,y,t)\in D\}=\{y\in\overline{\Omega}\,:\,(x,y,t)\in D_1\cup D_2\cup D_4\},\quad t\ge 1, 
\end{equation}
it follows from Theorem~\ref{Theorem:1.2} that 
\begin{equation}
\label{eq:2.20}
	\begin{aligned}
 	& 
	[\G(t)\phi](x) \ge C\int_\Omega \underline{H}(x,y,t)\phi^i(y)\,\dee y+C\int_{\partial\Omega} \underline{H}(x,y,t)\phi^b(y)\,\dee \sigma(y)
	\\
 	& 
	=C\int_{\{y\in\Omega \,:\,(x,y,t)\in D_1\}} P_N(x-y^*,t)\phi^i(y)\,\dee y\\
	 & \qquad
	 +C\int_{\{y\in\Omega \,:\,(x,y,t)\in (D_2\cup D_4)\}} \Gamma_N\left(x-y^*,\frac{t}{2}\right)\phi^i(y)\,\dee y\\
	& \qquad
	+C\int_{\{y\in\partial\Omega\,:\,(x,y,t)\in D_1\}} P_N(x-y^*,t)\phi^b(y)\,\dee\sigma(y)\\
	 & \qquad
	 +C\int_{\{y\in\partial\Omega\,:\,(x,y,t)\in (D_2\cup D_4)\}} \Gamma_N\left(x-y^*,\frac{t}{2}\right)\phi^b(y)\,\dee\sigma(y)
	\end{aligned}
\end{equation}
for $t\ge 1$.
If $(x,y,t)\in D_1$ and $t\ge 1$, then 
$$
	P_N(x-y^*,t)=c_N(x_N+y_N+t)|x-y^*+te_N|^{-N}\ge C
	\ge C\Gamma_N\left(x-y^*,\frac{t}{2}\right).
$$
This together with \eqref{eq:2.19} and \eqref{eq:2.20} implies that
\begin{equation*}
\begin{split}
	[\G(t)\phi](x)
	 & \ge C\int_{\{y\in\Omega\,:\,(x,y,t)\in (D_1\cup D_2\cup D_4)\}}\Gamma_N\left(x-y^*,\frac{t}{2}\right)
	\phi^i(y)\,\dee y\\
	 & \qquad
	 +C\int_{\{y\in\partial\Omega\,:\,(x,y,t)\in (D_1\cup D_2\cup D_4)\}}\Gamma_N\left(x-y^*,\frac{t}{2}\right)
	\phi^b(y)\,\dee\sigma(y)
	\\
 	& 
	=C\int_\Omega \Gamma_N\left(x-y^*,\frac{t}{2}\right)\phi^i(y)\,\dee y
	+C\int_{\partial\Omega} \Gamma_N\left(x-y^*,\frac{t}{2}\right)\phi^b(y)\,\dee \sigma(y).
\end{split}
\end{equation*}
Thus assertion~(2) holds, and the proof of Corollary~\ref{Corollary:1.1} is complete. 
$\Box$\vspace{5pt}
\newline
{\bf Proof of Corollary~\ref{Corollary:1.2}.}
We show assertion~(1).
Let $y\in\overline{\Omega}$ and $s\in(0,\infty)$. 
Set
$$
	u(x,t)
	:= G(x,y,t+s)-\int_\Omega G(x,z,t)G(z,y,s)\,\dee z-\int_{\partial\Omega}G(x,z,t)G(z,y,s)\,\dee\sigma(z)
$$
for $(x,t)\in\overline{\Omega}\times[0,\infty)$. 
By Theorem~\ref{Theorem:1.2} and Corollary~\ref{Corollary:1.1} 
we see that $u$ is a bounded and continuous solution to problem~\eqref{eq:H} with $\phi=(0,0)$. 
Then, by the comparison principle to bounded continuous solutions to problem~\eqref{eq:H} 
(see e.g., \cite{GH}*{Theorem~3.1}) we see that $u\equiv 0$ in $\overline{\Omega}\times[0,\infty)$. 
This implies the desired conclusion, and assertion~(1) follows.

We show assertion~(2).
Let $\phi^i$ and $\phi^b$ be nonnegative measurable functions in $\Omega$ and $\partial\Omega$, respectively, 
or $(\phi^i,\phi^b)\in L^q(\Omega)\times L^q(\partial\Omega)$ for some $q\in[1,\infty]$. 
Set $\phi=(\phi^i,\phi^b)$. 
Thanks to Corollary~\ref{Corollary:1.1}, 
it follows from Fubini's theorem and assertion~(1) that 
\begin{align*}
	& 
	[\G(t)\G(s)\phi](x)
	\\
 	& 
	=\int_\Omega G(x,z,t)\left(\int_\Omega G(z,y,s)\phi^i(y)\,\dee y
	+\int_{\partial\Omega}G(z,y,s)\phi^b(y)\,\dee\sigma(y)\right)\,\dee z
	\\
 	& \quad
  	+\int_{\partial\Omega} G(x,z,t)\left(\int_\Omega G(z,y,s)\phi^i(y)\,\dee y
	+\int_{\partial\Omega}G(z,y,s)\phi^b(y)\,\dee\sigma(y)\right)\,\dee\sigma(z)
	\\
 	& 
	=\int_\Omega\left(\int_\Omega G(x,z,t)G(z,y,s)\,\dee z
	+\int_{\partial\Omega}G(x,z,t)G(z,y,s)\,\dee\sigma(z)\right)\phi^i(y)\,\dee y
	\\
 	& \quad
 	+\int_{\partial\Omega}\left(\int_\Omega G(x,z,t)G(z,y,s)\,\dee z
	+\int_{\partial\Omega}G(x,z,t)G(z,y,s)\,\dee\sigma(z)\right)\phi^b(y)\,\dee \sigma(y)
	\\
 	& 
	=\int_\Omega G(x,y,t+s)\phi^i(y)\,\dee y+\int_{\partial\Omega}G(x,y,t+s)\phi^b(y)\,\dee\sigma(y)
	\\
 	& 
	=[\G(t+s)\phi](x),\quad x\in\overline{\Omega},\,\,t,s\in(0,\infty).
\end{align*}
Thus assertion~(2) follows, and the proof of Corollary~\ref{Corollary:1.2} is complete.
$\Box$\vspace{5pt}
\newline
{\bf Proof of Corollary~\ref{Corollary:1.3}.}
It follows from Corollary~\ref{Corollary:1.1}-(2) that 
\begin{equation}
\label{eq:2.21}
	\|\G(t)\|_{q\to r}
	\le\left\{
	\begin{array}{ll} 
	Ct^{-(N-1)\left(\frac{1}{q}-\frac{1}{r}\right)} & \mbox{for}\quad t\in(0,1],\vspace{3pt}\\
	Ct^{-\frac{N}{2}\left(\frac{1}{q}-\frac{1}{r}\right)} & \mbox{for}\quad t\in(1,\infty),
\end{array}
\right.
\end{equation}
where $1\le q\le r\le\infty$. 

Let $y=(y',y_N)\in\Omega$ with $y_N<1$, and set $g(x,t):=G(x,y,t)$ for $x\in\overline{\Omega}\times(0,\infty)$. 
Then it follows from Theorem~\ref{Theorem:1.1}-(3), 
Corollary~\ref{Corollary:1.1}, and Corollary~\ref{Corollary:1.2} that 
$$
\begin{aligned}
 	& 
	\|g(t)\|_{L^q(\Omega)\times L^q(\partial\Omega)}
	\\
 	& 
	=\|G(t/2)(g(t/2),g(t/2)|_{\partial\Omega})\|_{L^q(\Omega)\times L^q(\partial\Omega)}
	\\
 	& 
	\le C\left(t^{-(N-1)\left(1-\frac{1}{q}\right)}
	+t^{-\frac{N}{2}\left(1-\frac{1}{q}\right)}\right)\|(g(t/2),g(t/2)|_{\partial\Omega})\|_{L^1(\Omega)\times L^1(\partial\Omega)}
	\\
 	& 
	=C\left(t^{-(N-1)\left(1-\frac{1}{q}\right)}+t^{-\frac{N}{2}\left(1-\frac{1}{q}\right)}\right)
\end{aligned}
$$
for $t>0$, where $1\le q\le\infty$. 
On the other hand, 
by Theorem~\ref{Theorem:1.2} we see that
$$
{\sup\limits_{\substack{y=(y',y_N)\in\Omega,\\ y_N<1}}}\|g(t)\|_{L^\infty(\Omega)}\ge C{\sup\limits_{\substack{y=(y',y_N)\in\Omega,\\ y_N<1}}}\|\underline{H}(\cdot,y,t)\|_{L^\infty(\Omega)}
\ge
\left\{
\begin{array}{ll}
Ct^{-(N-1)} & \mbox{for}\quad t\in(0,1],\vspace{5pt}\\
Ct^{-\frac{N}{2}} & \mbox{for}\quad t\in(1,\infty).
\end{array}
\right.
%
%	\|g(t)\|_{L^\infty(\Omega)}\ge C\|\underline{H}(\cdot,y,t)\|_{L^\infty(\Omega)}
%	\ge
%	\left\{
%	\begin{array}{ll}
%	Ct^{-(N-1)} & \mbox{for}\quad t\in(0,1],\vspace{5pt}\\
%	Ct^{-\frac{N}{2}} & \mbox{for}\quad t\in(1,\infty).
%	\end{array}
%	\right.
$$
These imply that
\begin{equation}
\label{eq:2.22}
\begin{split}
\|\G(t)\|_{q\to \infty} & \ge {\sup\limits_{\substack{y=(y',y_N)\in\Omega,\\ y_N<1}}}\frac{\|\G(t)g(t)\|_{L^\infty(\Omega)\times L^\infty(\partial\Omega)}}{\|g(t)\|_{L^q(\Omega)\times L^q(\partial\Omega)}}
{=\sup\limits_{\substack{y=(y',y_N)\in\Omega,\\ y_N<1}}}
\frac{\|g(2t)\|_{L^\infty(\Omega)\times L^\infty(\partial\Omega)}}{\|g(t)\|_{L^q(\Omega)\times L^q(\partial\Omega)}}\\
 & \ge\left\{
\begin{array}{ll} 
Ct^{-\frac{N-1}{q}} & \mbox{for}\quad t\in(0,1],\vspace{3pt}\\
Ct^{-\frac{N}{2q}} & \mbox{for}\quad t\in(1,\infty),
\end{array}
\right.
\end{split}
%\begin{aligned}
%	\|\G(t)\|_{q\to \infty} 
%	& 
%	\ge \frac{\|\G(t)g(t)\|_{L^\infty(\Omega)\times L^\infty(\partial\Omega)}}{\|g(t)\|_{L^q(\Omega)\times L^q(\partial\Omega)}}
%	=\frac{\|g(2t)\|_{L^\infty(\Omega)\times L^\infty(\partial\Omega)}}{\|g(t)\|_{L^q(\Omega)\times L^q(\partial\Omega)}}
%	\\
% 	& 
%	\ge\left\{
%	\begin{array}{ll} 
%	Ct^{-\frac{N-1}{q}} & \mbox{for}\quad t\in(0,1],\vspace{3pt}\\
%	Ct^{-\frac{N}{2q}} & \mbox{for}\quad t\in(1,\infty),
%	\end{array}
%\right.
%\end{aligned}
\end{equation}
where $1\le q\le\infty$. 
Furthermore, by Corollary~\ref{Corollary:1.2}-(2) we have
\begin{equation}
\label{eq:2.23}
\|\G(t)\|_{q\to\infty}=\|\G(t/2)\G(t/2)\|_{q\to\infty}\le\|\G(t/2)\|_{q\to r} \|\G(t/2)\|_{r\to\infty},\quad t>0,
\end{equation}
where $1\le q\le r\le\infty$. 
Therefore,
by \eqref{eq:2.21} with $q=r$ and $r=\infty$, \eqref{eq:2.22}, and \eqref{eq:2.23} 
we obtain the desired relation, and Corollary~\ref{Corollary:1.3} follows.
$\Box$ 
%%%%%%%%%%%%%%%%%%%%%%%%%%%%%%%%%%%%%%
%%%%%%%%%%%%%%%%%%%%%%%%%%%%%%%%%%%%%%
\section{Global-in-time positive solutions to problem~\eqref{eq:SH}}
%%%%%%%%%%%%%%%%%%%%%%%%%%%%%%%%%%%%%%
%%%%%%%%%%%%%%%%%%%%%%%%%%%%%%%%%%%%%%
In this section, 
as an application of our decay estimates of solutions to problem~\eqref{eq:H}, 
we study the existence/nonexistence of global-in-time solutions to problem~\eqref{eq:SH}. 
Before stating our result to problem~\eqref{eq:SH}, we recall the celebrated result by Fujita~\cite{F} on the existence/nonexistence 
of global-in-time positive solutions to the problem
\begin{equation}
\tag{F}
\label{eq:F}
\partial_t u-\Delta u=u^p,\quad x\in{\mathbb R}^N,\,\,\,t>0,
\qquad
u(x,0)=\phi(x),\quad x\in{\mathbb R}^N,
\end{equation}
where $N\ge 1$, $p>1$, and $\phi$ is a nonnegative measurable function in ${\mathbb R}^N$. 
In \cite{F}, Fujita proved: 
\begin{itemize}
  \item[(F1)] 
  if $1<p<p_F:=1+2/N$, problem~\eqref{eq:F} possesses no global-in-time positive solutions;
  \item[(F2)] 
  if $p>p_F$, problem~\eqref{eq:F} possesses a global-in-time positive solution for some nonnegative initial data $\phi$.
\end{itemize}
We call the critical exponent $p_F$ the Fujita exponent for problem~\eqref{eq:F}. 
Assertion~(F1) also holds for the critical case $p=p_F$. This was proved in \cite{Haya} for $N=1,2$ and in \cites{KST, S} for $N\ge 3$. 
(See e.g., \cites{AW, W} for alternative proofs.)
Subsequently, the Fujita exponent has been identified extensively for various nonlinear parabolic problems by many researchers 
(see e.g., two surveys~\cites{DL, L}). 

Although the semilinear heat equation with the dynamical boundary condition has been studied in several papers 
(see e.g., \cites{BBR, BC, BP01, BP02, PM, R}), 
the Fujita exponent for problem~\eqref{eq:SH} is still open even in the case of $\Omega={\mathbb R}^N_+$.  
In this section we identify the Fujita exponent for problem~\eqref{eq:SH}. 
\medskip

We formulate the definition of solutions to problem~\eqref{eq:SH}.
\begin{definition}
\label{Definition:3.1}
	Let $\phi^i$ and $\phi^b$ be nonnegative measurable functions in $\Omega$ and $\partial\Omega$, respectively. 
	Let $T\in(0,\infty]$ and set $\phi:=(\phi^i,\phi^b)$.
	For any nonnegative continuous function~$u$ in $\overline{\Omega}\times(0,T)$, 
	we call $u$ a solution of problem~\eqref{eq:SH} in $\Omega\times(0,T)$
	if $u\in L^\infty_{{\rm loc}}((0,T):L^\infty({\mathbb R}^N))$ and $u$ satisfies
	\begin{equation}
	\label{eq:3.1}
		u(x,t)
		=[\G(t)\phi](x)+\int_0^t[\G(t-s)(u(s)^p,0)](x)\,\dee s,
		\quad (x,t)\in\overline{\Omega}\times(0,T).
	\end{equation}
In the case of $T=\infty$, we call $u$ a global-in-time solution to problem~\eqref{eq:SH}.
\end{definition} 
\begin{remark}
\label{Remark:3.1}
{\rm (1)} Let $\phi^i\in L^\infty(\Omega)\cap C(\overline{\Omega})$ be such that 
$\phi^i$ coincides with $\phi^b$ on $\partial\Omega$. 
Let $u$ be a classical bounded solution to problem~\eqref{eq:SH} in $\Omega\times(0,T)$ 
such that $u\in C(\overline{\Omega}\times[0,T))$ and $u(x,0)=\phi^i(x)$ on $\overline{\Omega}$. 
Then, for any $(x,t)\in\overline\Omega\times(0,T)$ and $s\in(0,t)$, 
since
\begin{align*}
 & \int_\Omega -(\Delta_y G)(x,y,t-s)u(y,s)\,\dee y+\int_\Omega G(x,y,t-s)(\Delta u)(y,s)\,\dee y\\
 & =-\int_{\partial\Omega}\langle\nabla_yG(x,y,t-s),\nu\rangle u(y,s)\,\dee\sigma(y)
 +\int_{\partial\Omega} G(x,y,t-s)\langle\nabla u(y,s),\nu\rangle\,\dee\sigma(y)\\
  & =\int_{\partial\Omega} (\partial_tG)(x,y,t-s) u(y,s)\,\dee\sigma(y)
 -\int_{\partial\Omega} G(x,y,t-s)(\partial_t u)(y,s)\,\dee\sigma(y)\\
 & =-\frac{\partial}{\partial s}\int_{\partial\Omega} G(x,y,t-s) u(y,s)\,\dee\sigma(y),
\end{align*}
we have
\begin{align*}
 	& 
	\frac{\partial}{\partial s}\int_\Omega G(x,y,t-s)u(y,s)\,\dee y
	\\
 	&
	 =\int_\Omega \{-(\partial_t G)(x,y,t-s)u(y,s)+G(x,y,t-s)(\partial_tu)(y,s)\}\,\dee y
	 \\
 	&
	 =\int_\Omega \{-(\Delta_y G)(x,y,t-s)u(y,s)+G(x,y,t-s)((\Delta u)(y,s)+u(y,s)^p)\}\,\dee y\\
	  & =-\frac{\partial}{\partial s}\int_{\partial\Omega} G(x,y,t-s) u(y,s)\,\dee\sigma(y)
	  +\int_\Omega G(x,y,t-s)u(y,s)^p\,\dee y. 
\end{align*}
Since $u(x,0)=\phi^i(x)$ for $x\in\Omega$ and $u(x,0)=\phi^b(x)$ for $x\in\partial\Omega$,
by Theorem~{\rm\ref{Theorem:1.1}-(1)} we obtain 
\begin{align*}
	u(x,t) 
	& 
	=\int_\Omega G(x,y,t)\phi^i(y)\,\dee y
	+\int_{\partial\Omega} G(x,y,t)\phi^b(y)\,\dee \sigma(y)
	+\int_0^t\int_\Omega G(x,y,t-s)u(y,s)^p\,\dee y\,\dee s
	\\
	& 
	=[\G(t)\phi](x)+\int_0^t[\G(t-s)(u(s)^p,0)](x)\,\dee s
\end{align*}
for $(x,t)\in\overline{\Omega}\times(0,T)$. 
Then $u$ satisfies \eqref{eq:3.1} for $(x,t)\in\overline\Omega\times(0,T)$, 
and the definition of solutions to problem~\eqref{eq:SH} in Definition~{\rm\ref{Definition:3.1}} is valid. 
\vspace{5pt}
\newline
{\rm (2)}
Let $\phi^i$ and $\phi^b$ be nonnegative measurable functions in~$\Omega$ and $\partial\Omega$, respectively. 
Let $u$ be a solution to problem~\eqref{eq:SH} in $\Omega\times(0,T)$ for some $T\in(0,\infty]$. 
For any $\tau\in(0,T)$, set 
$$
	v(x,t):=u(x,t+\tau)\mbox{ in $\overline{\Omega}\times[0,T-\tau)$},
	\qquad \phi^i_\tau:=u(\tau),\qquad \phi^b_\tau:=u(\tau)|_{\partial\Omega}. 
$$
Then $v$ is a solution to problem~\eqref{eq:SH} in $\Omega\times(0,T-\tau)$ with $\phi^i$ and $\phi^b$ replaced by $\phi^i_\tau$ and $\phi^b_\tau$, respectively.
Indeed, by Corollary~{\rm \ref{Corollary:1.2}-(2)} we have 
\begin{align*}
 	&
	\G(t)\phi_\tau+\int_0^t \G(t-s)(v(s)^p,0)\,\dee s
	\\
 	& 
	=\G(t)\left[\G(\tau)\phi+\int_0^\tau \G(\tau-s)(u(s)^p,0)\,\dee s\right]
	+\int_0^t\G(t-s)(v(s)^p,0)\,\dee s
	\\
 	& 
	=\G(t+\tau)\phi+\int_0^\tau \G(t+\tau-s)(u(s)^p,0)\,\dee s+\int_0^t\G(t-s)(u(s+\tau)^p,0)\,\dee s
	\\
 	& 
	=\G(t+\tau)\phi+\int_0^{t+\tau} \G(t+\tau-s)(u(s)^p,0)\,\dee s
	\\
 	& 
	=u(t+\tau)=v(t)
\end{align*}
for $t\in(0,T-\tau)$.
\vspace{5pt}
\newline
{\rm (3)} 
Let $u$ be a solution to problem~\eqref{eq:SH} in $\Omega\times(0,T)$, where $T\in(0,\infty]$. 
Since $u\in L^\infty_{{\rm loc}}((0,T):L^\infty({\mathbb R}^N))$, 
it follows from the smoothness of the kernel $G$, 
the decay of the kernel $\G$ {\rm ({\it see Theorem}~{\rm\ref{Theorem:1.2}} {\it and} \eqref{eq:1.1})}, 
and Remark~{\rm\ref{Remark:3.1}-(2)} that 
$u\in C^{2;1}(\overline{\Omega}\times(0,T))$ and $u$ satisfies 
$$
	\left\{
	\begin{array}{ll}
	\displaystyle{\partial_t u=\Delta u+u^p}, & (x,t)\in\Omega\times(0,T),\vspace{5pt}\\
	\displaystyle{\partial_tu+\partial_\nu u=0}, & (x,t)\in\partial\Omega\times(0,T),
	\end{array}
\right.
$$
pointwisely. 
\end{remark}
In this section we prove the following result on the Fujita exponent 
for problem~\eqref{eq:SH}.
\begin{theorem}
\label{Theorem:3.1}
	Let $N\ge 1$ and $p>1$. 
	Let $\phi^i$ and $\phi^b$ be nonnegative measurable functions in~$\Omega$ and $\partial\Omega$, respectively. 
	\begin{itemize}
	\item[\rm(1)]
	If $1<p\le p_F$, then problem~\eqref{eq:SH} possesses no global-in-time positive solutions.
	\item[\rm(2)]
	Let $p>p_F$.
	Then there exists $\delta>0$ such that, 
	if $\phi=(\phi^i,\phi^b)$ satisfies 
	$$
	\|\phi\|_{L^1(\Omega)\times L^1(\partial\Omega)}
	+\|\phi\|_{L^\infty(\Omega)\times L^\infty(\partial\Omega)}
	\le\delta,
	$$ 
	then problem~\eqref{eq:SH} possesses a global-in-time positive solution.
	\end{itemize}
\end{theorem}
%
%%%%%%%%%%%
\subsection{Proof of Theorem~\ref{Theorem:3.1}-(1)}
%%%%%%%%%%%
For the proof of Theorem~\ref{Theorem:3.1}-(1), we introduce some notation and prepare two lemmas.
Let $\mu$ be the first Dirichlet eigenvalue for $-\Delta$ in the annulus $E:=\{x\in{\mathbb R}^N\,:\,1<|x|<2\}$ 
and $\psi$ the corresponding eigenfunction such that $\psi>0$ in $E$ and $\|\psi\|_{L^1(\Omega)}=1$. 
For any $n\in\{1,2,\dots\}$, set 
\begin{align}
\label{eq:3.2}
 & E_n:=3ne_N+nE=\{x\in{\mathbb R}^N\,:\,n<|x-3ne_N|<2n\}\subset\Omega,\\
 \nonumber
 & \psi_n(x):=n^{-N}\psi(n^{-1}(x-3e_N))\quad\mbox{for}\quad x\in E_n. 
\end{align}
Then $\psi_n$ satisfies 
\begin{equation}
\label{eq:3.3}
	-\Delta\psi_n=n^{-2}\mu\psi_n\,\,\,\mbox{in}\,\,\,E_n,
	\quad
	\psi_n>0\,\,\,\mbox{in}\,\,\,E_n,
	\quad
	\psi_n=0\,\,\,\mbox{on}\,\,\,\partial E_n,
	\quad
	\|\psi_n\|_{L^1(E_n)}=1,
\end{equation}
for $n=1,2,\dots$. 
\begin{lemma}
\label{Lemma:3.1}
	Assume that problem~\eqref{eq:SH} possesses a global-in-time positive solution~$u$. 
	Then there exists $C_*>0$ such that
	\begin{equation}
	\label{eq:3.4}
		\left(\int_{E_n}u(x,2n^2+\tau)\psi_n(x)\,\dee x\right)^{p-1}\le C_*n^{-2}
	\end{equation}
	for $\tau\ge 0$ and $n=1,2,\dots$. 
\end{lemma}
{\bf Proof.}
Let $u$ be a global-in-time positive solution to problem~\eqref{eq:SH}. 
For any $\tau\ge 0$ and $n=1,2,\dots$, set 
$$
	w_n(x,t):=u(x,t+2n^2+\tau),\quad (x,t)\in\Omega\times[0,\infty).
$$
Thanks to Remark~\ref{Remark:3.1}-(2), (3), 
we see that 
$w_n\in C^{2;1}(\overline{\Omega}\times[0,\infty))$ and it 
is a classical positive solution to problem~\eqref{eq:SH} with $w_n(0)=u(\cdot,2n^2+\tau)$. 
Then, for any $n=1,2,\dots$, the function  
$$
	Z_n(t):=\int_{E_n}w_n(x,t)\psi_n(x)\,\dee x
$$ 
can be defined for $t\ge 0$ and it is positive in $[0,\infty)$.  
Since $\partial_\nu\psi_n<0$ on $\partial E_n$ by Hopf's lemma,  
it follows form \eqref{eq:3.3} and the Jensen inequality that
\begin{equation}
\label{eq:3.5}
	\begin{aligned}
	\frac{\dee}{\dee t}Z_n(t)
 	& 
	=\int_{E_n}\bigg(\Delta w_n+w_n^p\bigg)\psi_n\,\dee x
	\\
 	& 
	\ge-\int_{\partial E_n}w_n\partial_\nu\psi_n\,\dee\sigma
	+\int_{E_n}w_n\Delta\psi_n\,\dee x+\left(\int_{E_n}w_n\psi_n\,\dee x\right)^p
	\\
 	& 
	\ge -n^{-2}\mu Z_n(t)+Z_n(t)^p,\quad t>0. 
	\end{aligned}
\end{equation}
Assume that 
\begin{equation}
\label{eq:3.6}
	Z_{n_*}(0)^{p-1}>2n_*^{-2}\mu 
	\quad\mbox{for some $n_*\in\{1,2,\dots\}$}.
\end{equation}
By the continuity of $Z_{n_*}$ in $[0,\infty)$
we see that
\begin{equation}
\label{eq:3.7}
	T:=\sup\left\{t>0\,:\,Z_{n_*}(\tau)^{p-1}>\mu n_*^{-2}\mbox{ for $\tau\in(0,t)$}\right\}>0. 
\end{equation}
Then, by \eqref{eq:3.5} we have 
\begin{equation}
\label{eq:3.8}
	\frac{\dee}{\dee t}Z_{n_*}(t)\ge 0,\quad t\in(0,T),
\end{equation}
which together with \eqref{eq:3.6} implies that 
$$
	Z_{n_*}(t)^{p-1}\ge Z_{n_*}(0)^{p-1}> 2\mu n_*^{-2},\quad t\in(0,T). 
$$
We observe from \eqref{eq:3.6}, \eqref{eq:3.7}, and \eqref{eq:3.8} that $T=\infty$ and 
\begin{equation}
\label{eq:3.9}
	Z_{n_*}(t)^{p-1}\ge Z_{n_*}(0)^{p-1}>2\mu n_*^{-2}, \quad t>0.
\end{equation}
Therefore, by \eqref{eq:3.5} and \eqref{eq:3.9} we obtain 
$$
	\frac{d}{dt}Z_{n_*}(t)>\frac{1}{2}Z_{n_*}(t)^p, \quad t>0,
$$
which implies that $Z_{n_*}$ blows up in finite time. 
This is a contradiction. Therefore we see that 
$Z_n(0)^{p-1}\le 2\mu n^{-2}$ for $n=1,2,\dots$. 
This implies \eqref{eq:3.4}, and Lemma~\ref{Lemma:3.1} follows.
$\Box$
\begin{lemma}
\label{Lemma:3.2}
	Let $\zeta\in C_0({\mathbb R}^N)\setminus\{0\}$ be such that $\zeta\ge 0$ in ${\mathbb R}^N$, and set 
	\begin{equation*}
	[S(t)\zeta](x):=\int_{\mathbb R^N}\Gamma_N(x-y,t)\zeta(y)\,\dee y,\quad
	(x,t)\in{\mathbb R}^N\times(0,\infty). 
	\end{equation*}
	Then, for any compact set $K$ of ${\mathbb R}^N$, 
	there exists $C>0$ such that 
	$$
	[S(t)\zeta](\sqrt{t}x)
	\ge Ct^{-\frac{N}{2}}\int_{{\mathbb R}^N}\zeta(y)\,\dee y
	$$
	for $(x,t)\in K\times(1,\infty)$. 
\end{lemma}
{\bf Proof.}
Let $K$ be a compact set of ${\mathbb R}^N$. 
Since 
$$
	\bigg|x-\frac{y}{\sqrt t}\bigg|^2
	\le \bigg(|x|+\frac{|y|}{\sqrt t}\bigg)^2
	\le (|x|+|y|)^2<\infty
$$
for $x\in K$, $y\in\mbox{supp}\,\zeta$, and $t\ge 1$,
it follows from the definition of $\Gamma_N$ (see \eqref{eq:1.1}) that 
$$
	\Gamma_N(\sqrt t x-y,t)\ge Ct^{-\frac{N}{2}}
$$
for $x\in K$, $y\in\mbox{supp}\,\zeta$, and $t\ge 1$.
Then Lemma~\ref{Lemma:3.2} easily follows. 
$\Box$
\vspace{5pt}

Now we are ready to prove Theorem~\ref{Theorem:3.1}-(1). 
\vspace{5pt}
\newline
\noindent
{\bf Proof of Theorem~\ref{Theorem:3.1}-(1).}
Let $u$ be a global-in-time positive solution to problem~\eqref{eq:SH}. 
Then, thanks to Remark~\ref{Remark:3.1}-(2), 
we can assume, without loss of generality, that $\phi^i>0$ in $\Omega$. 
For $x=(x',x_N)\in \mathbb R^N$, set
$$
	\Phi(x)=
	\left\{
	\begin{array}{ll}
	0
	&\mbox{if}\quad x_N\ge0,
	\vspace{5pt}
	\\
	\phi^i(x',-x_N)
	&\mbox{if}\quad x_N<0.
	\end{array}
	\right.
$$
Then, by Corollary~\ref{Corollary:1.1}-(3) and \eqref{eq:3.1} we have
\begin{equation*}
	\begin{aligned}
	u(x,2t+2)
	&
	\ge [\G(2t+2)(\phi^i,0)](x)
	\\
	&
	\ge C \int_{\Omega}\Gamma_N\left(x-y^*, t+1\right)\phi^i(y)\,\dee y
	= C \int_{{\mathbb R}^N}\Gamma_N\left(x-y, t+1\right)\Phi(y)\,\dee y
	\end{aligned}
\end{equation*}
for $x\in \Omega$ and $t>0$. 
Thanks to the relation 
$$
\int_{{\mathbb R}^N}\Gamma_N(x-z,t)\Gamma_N(z-y,s)\,\dee z=\Gamma_N(x-y,t+s),
\quad x,y\in{\mathbb R}^N,\,\,t,s>0,
$$
it follows from Fubini's theorem that
\begin{equation}
\label{eq:3.10}
u(x,2t+2)\ge C\int_{{\mathbb R}^N}\Gamma_N(x-z,t)
	\left(\int_{{\mathbb R}^N}\Gamma_N(z-y,1)\Phi(y)\,\dee y\right)\,\dee z
\end{equation}
for $x\in \Omega$ and $t>0$. 
Furthermore, 
we find $\psi\in C_0({\mathbb R}^N)\setminus\{0\}$ such that 
$$
	\int_{{\mathbb R}^N}\Gamma_N(z-y,1)\Phi(y)\,\dee y\ge\psi(x)\ge 0,\quad x\in{\mathbb R}^N.
$$
Then, by Lemma~\ref{Lemma:3.2} and \eqref{eq:3.10}, 
for any compact set $K$ of $\Omega$, we have 
\begin{equation}
\label{eq:3.11}
u(x,2t+2)\ge C[S(t)\psi](x)\ge CMt^{-\frac{N}{2}}
\end{equation}
for $(x,t)\in\Omega\times(1,\infty)$ with $x\in \sqrt{t}K$, where $M:=\int_{{\mathbb R}^N}\psi(x)\,\dee x>0$. 

Consider the case of $1<p<p_F=1+2/N$. 
It follows from \eqref{eq:3.2} and \eqref{eq:3.11} that
$$
	u(x, 2n^2+2) \ge CMn^{-N},\quad x\in E_n,\,\, n=1,2,\dots.
$$
This together with \eqref{eq:3.3} implies that
$$
	\int_{E_n}u(x,2n^2+2)\psi_n(x)\,\dee x
	\ge CMn^{-N},\quad n=1,2,\dots. 
$$
Since $p<p_F$, it turns out that
$$
	n^2\left(\int_{E_n}u(x,2n^2+2)\psi_n(x)\,\dee x\right)^{p-1}\to\infty
	\quad\mbox{as}\quad n\to\infty,
$$
which contradicts \eqref{eq:3.4}. 
This means that problem~\eqref{eq:SH} possesses no global-in-time positive solutions. 
Thus Theorem~\ref{Theorem:3.1}-(1) follows in case of $1<p<p_F$. 

Consider the case of $p=p_F$. 
It follows that
\begin{equation}
\label{eq:3.12}
	\Gamma_N\left(x-y^*,\frac{t-s}{2}\right)
	=(2\pi(t-s))^{-\frac{N}{2}}e^{-\frac{|x-y^*|^2}{2(t-s)}}
	\ge Ct^{-\frac{N}{2}}
\end{equation}
for $x$, $y\in\Omega$ and $0<s<t$ with
$$
	\sqrt{t}\le|x-3\sqrt{t}e_N|\le 2\sqrt{t},\quad
	\sqrt{s}\le|y-3\sqrt{s}e_N|\le 2\sqrt{s},\quad
	1\le s\le \frac{t}{2},\qquad t>2.
$$
Since $G\ge H$ in $D$, 
by the definition of solutions to problem~\eqref{eq:SH} (see \eqref{eq:3.1}) we have 
\begin{align*}
u(x,2t+2) & 
\ge\int_4^{t+2}\int_\Omega G(x,y,2t+2-s)u(y,s)^p\,\dee y\,\dee s\\
 & \ge\int_4^{t+2}\int_\Omega H(x,y,2t+2-s)u(y,s)^p\,\dee y\,\dee s\\
 & \ge 2\int_1^{t/2}\int_{\{\sqrt{s}\le |y-3\sqrt{s}e_N|\le 2\sqrt{s}\}}H(x,y,2t-2s)u(y,2s+2)^p\,\dee y\,\dee s
\end{align*}
for $(x,t)\in\Omega\times(2,\infty)$. 
This together with Theorem~\ref{Theorem:1.2} in the case of $(x,y,t)\in D_2\cup D_4$ and \eqref{eq:3.12} 
implies that
$$
	\begin{aligned}
	u(x,2t+2)
	& 
	\ge
	C\int_1^{t/2}\int_{\{\sqrt{s}\le |y-3\sqrt{s}e_N|\le 2\sqrt{s}\}}
	\underline{H}(x-y^*,t-s)u(y,2s+2)^p\,\dee y\,\dee s
	\\
 	& 
	\ge
	C\int_1^{t/2}\int_{\{\sqrt{s}\le |y-3\sqrt{s}e_N|\le 2\sqrt{s}\}}
	\Gamma_N\left(x-y^*,\frac{t-s}{2}\right)u(y,2s+2)^p\,\dee y\,\dee s
	\\
 	& 
	\ge Ct^{-\frac{N}{2}}\int_1^{t/2} s^{-\frac{N}{2}(p-1)}\,\dee s\ge Ct^{-\frac{N}{2}}\log t
	\end{aligned}
$$
for $x\in\Omega$ and $t>8(N+2)$ with $\sqrt{t}\le|x-3\sqrt{t}e_N|\le 2\sqrt{t}$.
Then 
$$
	n^2\left(\int_{E_n}u(x,2n^2+2)\psi_n(x)\,\dee x\right)^{p-1}
	\ge n^2\left(Cn^{-N}\log n\right)^{\frac{2}{N}}
	\ge C(\log n)^{\frac{2}{N}}\to\infty
$$
as $n\to\infty$, which contradicts \eqref{eq:3.4}. 
Therefore, similarly to the case of $p<p_F$, we see that problem~\eqref{eq:SH} possesses no global-in-time positive solutions. 
Thus Theorem~\ref{Theorem:3.1}-(1) follows in case $p=p_F$,
and the proof of Theorem~\ref{Theorem:3.1}-(1) is complete.
$\Box$
%%%%%%%%%%%
\subsection{Proof of Theorem~\ref{Theorem:3.1}-(2)}
%%%%%%%%%%%
We apply the contraction mapping theorem to prove Theorem~\ref{Theorem:3.1}-(2).
\vspace{5pt}
\newline
{\bf Proof of Theorem~\ref{Theorem:3.1}-(2).}
Let $p>p_F$ and $\delta>0$ be small enough. 
Let 
$$
	\phi=(\phi^i,\phi^b)
	\in (L^1(\Omega)\times L^1(\partial\Omega))\cap (L^\infty(\Omega)\times L^\infty(\partial\Omega))
$$ 
be such that
$$
	\|\phi\|_{L^1(\Omega)\times L^1(\partial\Omega)}
	+\|\phi\|_{L^\infty(\Omega)\times L^\infty(\partial\Omega)}
	\le\delta.
$$ 
By Corollary~\ref{Corollary:1.1}-(2) we find $C_*>0$ such that
\begin{equation}
\label{eq:3.13}
	\begin{aligned}
	\|\G(t)\phi\|_{L^q(\Omega)} & \le C_*(1+t)^{-\frac{N}{2}\left(1-\frac{1}{q}\right)}
	\left(\|\phi\|_{L^1(\Omega)\times L^1(\partial\Omega)}
	+\|\phi\|_{L^\infty(\Omega)\times L^\infty(\partial\Omega)}\right)
	\\
 	& 
	\le C_*\delta(1+t)^{-\frac{N}{2}\left(1-\frac{1}{q}\right)}
	\end{aligned}
\end{equation}
for $t>0$ and $q\in[1,\infty]$.
{Define
$$
X:=\{f\in L^\infty((0,\infty), L^1(\Omega)\cap L^\infty(\Omega))\,:\,\|f\|_X<\infty\},
$$
where
$$
\|v\|_X:=\sup_{t>0}\,\|v(t)\|_{L^1(\Omega)}+\sup_{t>0}\, (1+t)^{\frac{N}{2}}\|v(t)\|_{L^\infty(\Omega)}.
$$}
Then $X$ is a Banach space equipped with the norm $\|\cdot\|_X$. 
Let $B$ be a closed ball in $X$ defined by  
$$
	B:=\{v\in X\,:\,\|v\|_X\le 4C_*\delta\}.
$$
Let $v_1$, $v_2\in B$. Since
$$
	|v_1(x,t)^p-v_2(x,t)^p|\le p(|v_1(x,t)|+|v_2(x,t)|)^{p-1}|v_1(x,t)-v_2(x,t)|
$$
for almost all $(x,t)\in \Omega\times(0,\infty)$, 
we have 
\begin{equation}
\label{eq:3.14}
	\begin{aligned}
 	& \|v_1(t)^p-v_2(t)^p\|_{L^r(\Omega)}
	\\
 	& 
	\le p\left(\|v_1(t)\|_{L^\infty(\Omega)}+\|v_2(t)\|_{L^\infty(\Omega)}\right)^{p-1}
	\|v_1(t)-v_2(t)\|_{L^r(\Omega)}
	\\
 	& 
	\le p(8C_*\delta)^{p-1}(1+t)^{-\frac{N}{2}(p-1)}
	\|v_1(t)-v_2(t)\|_{L^1(\Omega)}^{\frac{1}{r}}\|v_1(t)-v_2(t)\|_{L^\infty(\Omega)}^{1-\frac{1}{r}}
	\\
	& 
	\le C(C_*\delta)^{p-1}(1+t)^{-\frac{N}{2}\left(p-\frac{1}{r}\right)}\|v_1-v_2\|_X
	\end{aligned}
\end{equation}
for $t>0$ and $r\in[1,\infty]$. 
Since $p>p_F=1+2/N$, 
it follows from \eqref{eq:1.8} and \eqref{eq:3.14} that
\begin{equation}
\label{eq:3.15}
	\begin{aligned}
	& 
	\left\|\,\int_0^t \G(t-s)(v_1(s)^p-v_2(s)^p,0)\,\dee s\,\right\|_{L^q(\Omega)}
	\\
	& 
	\le\int_0^t\|\G(t-s)(v_1(s)^p-v_2(s)^p,0)\|_{L^q(\Omega)}\,\dee s
	\le \int_0^t\|v_1(s)^p-v_2(s)^p\|_{L^{q}(\Omega)}\,\dee s
	\\
	& 
	\le C(C_*\delta)^{p-1}\|v_1-v_2\|_X\int_0^\infty(1+s)^{-\frac{N}{2}\left(p-\frac{1}{q}\right)}\,\dee s
	\le C(C_*\delta)^{p-1}\|v_1-v_2\|_X
	\end{aligned}
\end{equation}
for $t>0$ and $q\in[1,\infty]$. 
Similarly,
by \eqref{eq:1.7}, \eqref{eq:1.8}, and \eqref{eq:3.14} we see that
\begin{equation}
\label{eq:3.16}
	\begin{aligned}
	& 
	\left\|\,\int_0^t \G(t-s)(v_1(s)^p-v_2(s)^p,0)\,\dee s\,\right\|_{L^q(\Omega)}
	\\
	& \le \left(\int_0^{t/2}+\int_{t/2}^t\right)\|\G(t-s)(v_1(s)^p-v_2(s)^p,0)\|_{L^q(\Omega)}\,\dee s
	\\
	&
	\le C\int_0^{t/2}(t-s)^{-\frac{N}{2}\left(1-\frac{1}{q}\right)}
	\|v_1(s)^p-v_2(s)^p\|_{L^1(\Omega)}\,\dee s
	\\
	&\hspace{5cm}
	+\int_{t/2}^t\|v_1(s)^p-v_2(s)^p\|_{L^q(\Omega)}\,\dee s
	\\
	&
	\le C(C_*\delta)^{p-1}\|v_1-v_2\|_X
	\\
	& \qquad\quad
	\times\biggr[\left(\frac{t}{2}\right)^{-\frac{N}{2}\left(1-\frac{1}{q}\right)}
	\int_0^\infty (1+s)^{-\frac{N}{2}(p-1)}\,\dee s
	+\int_{t/2}^\infty (1+s)^{-\frac{N}{2}\left(p-\frac{1}{q}\right)}\,\dee s\biggr]
	\\
	&
	\le C(C_*\delta)^{p-1} \|v_1-v_2\|_X (1+t)^{-\frac{N}{2}\left(1-\frac{1}{q}\right)}
	\end{aligned}
\end{equation}
for $t\ge 2$ and $q\in[1,\infty]$.
Here we used the relation that $t-s\ge 1$ for $s\in(0,t/2)$ and $t\ge 2$.
Combining \eqref{eq:3.15} and \eqref{eq:3.16}, we obtain 
\begin{equation}
\label{eq:3.17}
	\left\|\int_0^t \G(t-s)(v_1(s)^p-v_2(s)^p,0)\,\dee s\,\right\|_{L^q(\Omega)}
	\le C(C_*\delta)^{p-1} \|v_1-v_2\|_X (1+t)^{-\frac{N}{2}\left(1-\frac{1}{q}\right)}
\end{equation}
for $t>0$ and $q\in[1,\infty]$. 

For $v\in B$, set
$$
	Q[v](t):=\G(t)\phi+\int_0^t \G(t-s)(v(s)^p,0)\,\dee s,\quad t>0.
$$
Then, 
taking small enough $\delta>0$ if necessary, 
by \eqref{eq:3.13} and \eqref{eq:3.17} with $v_1=v$ and $v_2=0$ we have
$$
	\begin{aligned}
	 & \|Q[v]\|_X \\
	& 
	\le\sup_{t>0}\,\|\G(t)\phi\|_{L^1(\Omega)}
	+\sup_{t>0}\,(1+t)^{\frac{N}{2}}\|\G(t)\phi\|_{L^\infty(\Omega)}
	\\
 	&\quad
	+\sup_{t>0}\,\left\|\,\int_0^t \G(t-s)(v(s)^p,0)\,\dee s\,\right\|_{L^1(\Omega)}
 	+\sup_{t>0}\,(1+t)^{\frac{N}{2}}\left\|\,\int_0^t \G(t-s)(v(s)^p,0)\,\dee s\,\right\|_{L^\infty(\Omega)}
	\\
  	& 
	\le 2C_*\delta+2C(C_*\delta)^{p-1}\|v\|_X\le 4C_*\delta,
	\end{aligned}
$$
which implies that $Q[v]\in B$ for $v\in B$.
Furthermore, taking small enough $\delta>0$ again if necessary, 
by \eqref{eq:3.17} we have
$$
	\begin{aligned}
	\left\|Q[v_1]-Q[v_2]\right\|_X
 	& 
	\le\sup_{t>0}\,\left\|\int_0^t \G(t-s)(v_1(s)^p-v_2(s)^p,0)\,\dee s\,\right\|_{L^1(\Omega)}
	\\
 	& \qquad\qquad
 	+\sup_{t>0}\,(1+t)^{\frac{N}{2}}
	\left\|\int_0^t \G(t-s)(v_1(s)^p-v_2(s)^p,0)\,\dee s\,\right\|_{L^\infty(\Omega)}
	\\
 	& 
	\le 2C(C_*\delta)^{p-1} \|v_1-v_2\|_X\le\frac{1}{2}\|v_1-v_2\|_X
	\end{aligned}
$$
for $v_1$, $v_2\in B$. 
Then we apply the contraction mapping theorem
to find a unique $u\in B$ such that
$$
	u(t)=Q[u](t)=\G(t)\phi+\int_0^t\G(t-s)(u(s)^p,0)\,\dee s
	\quad\mbox{in}\quad X.
$$
Then, since $u\in L^\infty(\Omega\times(0,\infty))$, 
thanks to the smoothness of the kernel $\G$ and the decay of $\G$ (see Theorem~\ref{Theorem:1.2} and \eqref{eq:1.1}), 
we see that $u\in C(\overline{\Omega}\times(0,\infty))$. 
Thus $u$ is a global-in-time solution to problem~\eqref{eq:SH}, and  Theorem~\ref{Theorem:3.1}-(2) follows. 
The proof of Theorem~\ref{Theorem:3.1} is complete. 
$\Box$
\vspace{8pt}

\noindent
{\bf Acknowledgment.}
T. K. was supported in part by JSPS KAKENHI Grant Number JP 22KK0035.
K. I. and T. K. were supported in part by JSPS KAKENHI Grant Number 19H05599. 
They would also like to thank Professor Marek Fila for intending our attention to 
elliptic and parabolic problems with the homogeneous dynamical boundary condition
and for his friendship during his lifetime.
%
%%%%%%%%%%%%%%%%%%%%%%%%%%%%%%%%%%%%%%
%%%%%%%%%%%%    references    %%%%%%%%%%%%%%%%%%
%%%%%%%%%%%%%%%%%%%%%%%%%%%%%%%%%%%%%%

%\begin{thebibliography}{10}
%
%\bibitem{E2}
%J. Escher, 
%Quasilinear parabolic systems with dynamical boundary conditions, 
%Commun. Partial Differ. Equ. 18 (1993), 1309--1364.
%%
%\bibitem{E3}
%J. Escher,
%On the qualitative behaviour of some semilinear parabolic problems,
%Differ. Integral Equ. 8 (1995), 247--267.
%%
%\bibitem{F}
%A. Friedman,  
%{\it Partial differential equations of parabolic type}, 
%Prentice-Hall, Inc., Englewood Cliffs, N.J. (1964).
%%
%\bibitem{GM}
%C. G. Gal and M. Meyries, 
%Nonlinear elliptic problems with dynamical boundary conditions of reactive and reactive-diffusive type,
%Proc. Lond. Math. Soc. 108 (2014), 1351--1380. 
%%
%\end{thebibliography}

\begin{bibdiv}
\begin{biblist}
%%%%%
\bib{AF}{article}{
   author={Amann, H.},
   author={Fila, M.},
   title={A Fujita-type theorem for the Laplace equation with a dynamical
   boundary condition},
   journal={Acta Math. Univ. Comenian. (N.S.)},
   volume={66},
   date={1997},
%   number={2},
   pages={321--328},
%   issn={0862-9544},
%   review={\MR{1620441}},
}
%%%%%
\bib{AQR}{article}{
   author={Arrieta, Jos\'e M.},
   author={Quittner, Pavol},
   author={Rodr\'iguez-Bernal, An\'ibal},
   title={Parabolic problems with nonlinear dynamical boundary conditions
   and singular initial data},
   journal={Differential Integral Equations},
   volume={14},
   date={2001},
%   number={12},
   pages={1487--1510},
%   issn={0893-4983},
%   review={\MR{1859918}},
}
%%%%%
\bib{AW}{article}{
   author={Aronson, D. G.},
   author={Weinberger, H. F.},
   title={Multidimensional nonlinear diffusion arising in population
   genetics},
   journal={Adv. in Math.},
   volume={30},
   date={1978},
%   number={1},
   pages={33--76},
%   issn={0001-8708},
%   review={\MR{0511740}},
%   doi={10.1016/0001-8708(78)90130-5},
}
%%%%%
\bib{BBR}{article}{
   author={Bandle, Catherine},
   author={von Below, Joachim},
   author={Reichel, Wolfgang},
   title={Parabolic problems with dynamical boundary conditions: eigenvalue
   expansions and blow up},
   journal={Atti Accad. Naz. Lincei Rend. Lincei Mat. Appl.},
   volume={17},
   date={2006},
%   number={1},
   pages={35--67},
%   issn={1120-6330},
%   review={\MR{2237743}},
%   doi={10.4171/RLM/453},
}
%%%%%
\bib{BC}{article}{
   author={von Below, Joachim},
   author={De Coster, Colette},
   title={A qualitative theory for parabolic problems under dynamical
   boundary conditions},
   journal={J. Inequal. Appl.},
   volume={5},
   date={2000},
%   number={5},
   pages={467--486},
%   issn={1025-5834},
%   review={\MR{1800975}},
%   doi={10.1155/S1025583400000266},
}
%%%%%
\bib{BP01}{article}{
   author={von Below, Joachim},
   author={Pincet Mailly, Ga\"elle},
   title={Blow up for reaction diffusion equations under dynamical boundary
   conditions},
   journal={Comm. Partial Differential Equations},
   volume={28},
   date={2003},
%   number={1-2},
   pages={223--247},
%   issn={0360-5302},
%   review={\MR{1974455}},
%   doi={10.1081/PDE-120019380},
}
%%%%%
\bib{BP02}{article}{
   author={von Below, Joachim},
   author={Pincet Mailly, Ga\"elle},
   title={Blow up for some nonlinear parabolic problems with convection
   under dynamical boundary conditions},
   journal={Discrete Contin. Dyn. Syst.},
   date={2007},
   pages={1031--1041},
%   issn={1078-0947},
%   isbn={978-1-60133-010-9; 1-60133-010-3},
%   review={\MR{2409941}},
}
%%%%%
\bib{C}{book}{
   author={Crank, J.},
   title={The mathematics of diffusion},
%   edition={2},
   publisher={Clarendon Press, Oxford},
   date={1975},
   pages={ix+414},
%   review={\MR{0359551}},
}
%%%%%
\bib{DL}{article}{
   author={Deng, Keng},
   author={Levine, Howard A.},
   title={The role of critical exponents in blow-up theorems: the sequel},
   journal={J. Math. Anal. Appl.},
   volume={243},
   date={2000},
%   number={1},
   pages={85--126},
%   issn={0022-247X},
%   review={\MR{1742850}},
%   doi={10.1006/jmaa.1999.6663},
}
%%%%%
\bib{DPZ}{article}{
   author={Denk, Robert},
   author={Pr\"uss, Jan},
   author={Zacher, Rico},
   title={Maximal $L_p$-regularity of parabolic problems with boundary
   dynamics of relaxation type},
   journal={J. Funct. Anal.},
   volume={255},
   date={2008},
%   number={11},
   pages={3149--3187},
%   issn={0022-1236},
%   review={\MR{2464573}},
%   doi={10.1016/j.jfa.2008.07.012},
}
%%%%%
\bib{EMR}{article}{
   author={ter Elst, A. F. M.},
   author={Meyries, M.},
   author={Rehberg, J.},
   title={Parabolic equations with dynamical boundary conditions and source
   terms on interfaces},
   journal={Ann. Mat. Pura Appl. (4)},
   volume={193},
   date={2014},
%   number={5},
   pages={1295--1318},
%   issn={0373-3114},
%   review={\MR{3262633}},
%   doi={10.1007/s10231-013-0329-7},
}
%%%%%
\bib{E}{article}{
   author={Escher, Joachim},
   title={Quasilinear parabolic systems with dynamical boundary conditions},
   journal={Comm. Partial Differential Equations},
   volume={18},
   date={1993},
%   number={7-8},
   pages={1309--1364},
%   issn={0360-5302},
%   review={\MR{1233197}},
%   doi={10.1080/03605309308820976},
}
%%%%%
\bib{FIK01}{article}{
   author={Fila, Marek},
   author={Ishige, Kazuhiro},
   author={Kawakami, Tatsuki},
   title={Convergence to the Poisson kernel for the Laplace equation with a
   nonlinear dynamical boundary condition},
   journal={Commun. Pure Appl. Anal.},
   volume={11},
   date={2012},
%   number={3},
   pages={1285--1301},
%   issn={1534-0392},
%   review={\MR{2968622}},
%   doi={10.3934/cpaa.2012.11.1285},
}
%%%%%
\bib{FIK02}{article}{
   author={Fila, Marek},
   author={Ishige, Kazuhiro},
   author={Kawakami, Tatsuki},
   title={Large-time behavior of solutions of a semilinear elliptic equation
   with a dynamical boundary condition},
   journal={Adv. Differential Equations},
   volume={18},
   date={2013},
%   number={1-2},
   pages={69--100},
%   issn={1079-9389},
%   review={\MR{3052711}},
}
%%%%%
\bib{FIK02a}{article}{
   author={Fila, Marek},
   author={Ishige, Kazuhiro},
   author={Kawakami, Tatsuki},
   title={Existence of positive solutions of a semilinear elliptic equation
   with a dynamical boundary condition},
   journal={Calc. Var. Partial Differential Equations},
   volume={54},
   date={2015},
%   number={2},
   pages={2059--2078},
%   issn={0944-2669},
%   review={\MR{3396444}},
%   doi={10.1007/s00526-015-0856-8},
}
%%%%%
\bib{FIK02b}{article}{
   author={Fila, Marek},
   author={Ishige, Kazuhiro},
   author={Kawakami, Tatsuki},
   title={Minimal solutions of a semilinear elliptic equation with a
   dynamical boundary condition},
%   language={English, with English and French summaries},
   journal={J. Math. Pures Appl.},
   volume={105},
   date={2016},
%   number={6},
   pages={788--809},
%   issn={0021-7824},
%   review={\MR{3491532}},
%   doi={10.1016/j.matpur.2015.11.014},
}
%%%%%
\bib{FIK03}{article}{
   author={Fila, Marek},
   author={Ishige, Kazuhiro},
   author={Kawakami, Tatsuki},
   title={The large diffusion limit for the heat equation with a dynamical
   boundary condition},
   journal={Commun. Contemp. Math.},
   volume={23},
   date={2021},
%   number={1},
   pages={Paper No. 2050003, 20},
%   issn={0219-1997},
%   review={\MR{4169250}},
%   doi={10.1142/s0219199720500030},
}
%%%%%
\bib{FIK04}{article}{
   author={Fila, Marek},
   author={Ishige, Kazuhiro},
   author={Kawakami, Tatsuki},
   title={Solvability of the heat equation on a half-space with a dynamical
   boundary condition and unbounded initial data},
   journal={Z. Angew. Math. Phys.},
   volume={74},
   date={2023},
%   number={4},
   pages={Paper No. 143, 17},
%   issn={0044-2275},
%   review={\MR{4605611}},
%   doi={10.1007/s00033-023-02040-7},
}
%%%%%
\bib{FIKL}{article}{
   author={Fila, Marek},
   author={Ishige, Kazuhiro},
   author={Kawakami, Tatsuki},
   author={Lankeit, Johannes},
   title={The large diffusion limit for the heat equation in the exterior of
   the unit ball with a dynamical boundary condition},
   journal={Discrete Contin. Dyn. Syst.},
   volume={40},
   date={2020},
%   number={11},
   pages={6529--6546},
%   issn={1078-0947},
%   review={\MR{4147360}},
%   doi={10.3934/dcds.2020289},
}
%%%%%
\bib{FQ}{article}{
   author={Fila, Marek},
   author={Quittner, Pavol},
   title={Large time behavior of solutions of a semilinear parabolic
   equation with a nonlinear dynamical boundary condition},
   conference={
      title={Topics in nonlinear analysis},
   },
   book={
      series={Progr. Nonlinear Differential Equations Appl.},
      volume={35},
      publisher={Birkh\"auser, Basel},
   },
%   isbn={3-7643-6016-X},
   date={1999},
   pages={251--272},
%   review={\MR{1725573}},
}
%%%%%
\bib{FV}{article}{
   author={Fiscella, Alessio},
   author={Vitillaro, Enzo},
   title={Local Hadamard well-posedness and blow-up for reaction-diffusion
   equations with non-linear dynamical boundary conditions},
   journal={Discrete Contin. Dyn. Syst.},
   volume={33},
   date={2013},
%   number={11-12},
   pages={5015--5047},
%   issn={1078-0947},
%   review={\MR{3060824}},
%   doi={10.3934/dcds.2013.33.5015},
}
%%%%%
%\bib{Fried}{book}{
%   author={Friedman, Avner},
%   title={Partial differential equations of parabolic type},
%   publisher={Prentice-Hall, Inc., Englewood Cliffs, NJ},
%   date={1964},
%   pages={xiv+347},
%%   review={\MR{0181836}},
%}
%%%%%
\bib{F}{article}{
   author={Fujita, Hiroshi},
   title={On the blowing up of solutions of the Cauchy problem for
   $u\sb{t}=\Delta u+u\sp{1+\alpha }$},
   journal={J. Fac. Sci. Univ. Tokyo Sect. I},
   volume={13},
   date={1966},
   pages={109--124 (1966)},
%   issn={0368-2269},
%   review={\MR{0214914}},
}
%%%%%
\bib{GM}{article}{
   author={Gal, Ciprian G.},
   author={Meyries, Martin},
   title={Nonlinear elliptic problems with dynamical boundary conditions of
   reactive and reactive-diffusive type},
   journal={Proc. Lond. Math. Soc. (3)},
   volume={108},
   date={2014},
%   number={6},
   pages={1351--1380},
%   issn={0024-6115},
%   review={\MR{3218312}},
%   doi={10.1112/plms/pdt057},
}
%%%%%
\bib{GH}{article}{
   author={Giga, Yoshikazu},
   author={Hamamuki, Nao},
   title={On a dynamic boundary condition for singular degenerate parabolic
   equations in a half space},
   journal={NoDEA Nonlinear Differential Equations Appl.},
   volume={25},
   date={2018},
%   number={6},
   pages={Paper No. 51, 39},
%   issn={1021-9722},
%   review={\MR{3863547}},
%   doi={10.1007/s00030-018-0542-6},
}
%%%%%
\bib{HQ}{article}{
   author={Hamamuki, Nao},
   author={Liu, Qing},
   title={A deterministic game interpretation for fully nonlinear parabolic
   equations with dynamic boundary conditions},
   journal={ESAIM Control Optim. Calc. Var.},
   volume={26},
   date={2020},
   pages={Paper No. 13, 42},
%   issn={1292-8119},
%   review={\MR{4064473}},
%   doi={10.1051/cocv/2019076},
}
%%%%%
\bib{Haya}{article}{
   author={Hayakawa, Kantaro},
   title={On nonexistence of global solutions of some semilinear parabolic
   differential equations},
   journal={Proc. Japan Acad.},
   volume={49},
   date={1973},
   pages={503--505},
%  issn={0021-4280},
%   review={\MR{0338569}},
}
%%%%%
\bib{H}{article}{
   author={Hintermann, Thomas},
   title={Evolution equations with dynamic boundary conditions},
   journal={Proc. Roy. Soc. Edinburgh Sect. A},
   volume={113},
   date={1989},
%   number={1-2},
   pages={43--60},
%   issn={0308-2105},
%   review={\MR{1025453}},
%   doi={10.1017/S0308210500023945},
}
%%%%%
\bib{IK}{article}{
   author={Ishige, Kazuhiro},
   author={Kawakami, Tatsuki},
   title={Critical Fujita exponents for semilinear heat equations with
   quadratically decaying potential},
   journal={Indiana Univ. Math. J.},
   volume={69},
   date={2020},
   number={6},
   pages={2171--2207},
%   issn={0022-2518},
%   review={\MR{4170090}},
%   doi={10.1512/iumj.2020.69.7989},
}
%%%%%
\bib{K}{article}{
   author={Kirane, M.},
   title={Blow-up for some equations with semilinear dynamical boundary
   conditions of parabolic and hyperbolic type},
   journal={Hokkaido Math. J.},
   volume={21},
   date={1992},
%   number={2},
   pages={221--229},
%   issn={0385-4035},
%   review={\MR{1169789}},
%   doi={10.14492/hokmj/1381413677},
}
%%%%%
\bib{KST}{article}{
   author={Kobayashi, Kusuo},
   author={Sirao, Tunekiti},
   author={Tanaka, Hiroshi},
   title={On the growing up problem for semilinear heat equations},
   journal={J. Math. Soc. Japan},
   volume={29},
   date={1977},
%   number={3},
   pages={407--424},
%   issn={0025-5645},
%   review={\MR{0450783}},
%   doi={10.2969/jmsj/02930407},
}
%%%%%
\bib{L}{article}{
   author={Levine, Howard A.},
   title={The role of critical exponents in blowup theorems},
   journal={SIAM Rev.},
   volume={32},
   date={1990},
%   number={2},
   pages={262--288},
%   issn={0036-1445},
%   review={\MR{1056055}},
%   doi={10.1137/1032046},
}
%%%%%
\bib{PM}{article}{
   author={Pincet Mailly, Ga\"elle},
   title={Blow up for nonlinear parabolic equations with time degeneracy
   under dynamical boundary conditions},
   journal={Nonlinear Anal.},
   volume={67},
   date={2007},
%   number={3},
   pages={657--667},
%   issn={0362-546X},
%   review={\MR{2319200}},
%   doi={10.1016/j.na.2006.06.021},
}
%%%%%
\bib{R}{article}{
   author={Rault, Jean-Fran\c cois},
   title={The Fujita phenomenon in exterior domains under dynamical boundary
   conditions},
   journal={Asymptot. Anal.},
   volume={66},
   date={2010},
%   number={1},
   pages={1--8},
%   issn={0921-7134},
%   review={\MR{2582445}},
}
%%%%%
\bib{S}{article}{
   author={Sugitani, Sadao},
   title={On nonexistence of global solutions for some nonlinear integral
   equations},
   journal={Osaka Math. J.},
   volume={12},
   date={1975},
   pages={45--51},
%   issn={0388-0699},
%   review={\MR{0470493}},
}
%%%%%
\bib{VV01}{article}{
   author={Vazquez, Juan Luis},
   author={Vitillaro, Enzo},
   title={Heat equation with dynamical boundary conditions of locally
   reactive type},
   journal={Semigroup Forum},
   volume={74},
   date={2007},
%   number={1},
   pages={1--40},
%   issn={0037-1912},
%   review={\MR{2301570}},
%   doi={10.1007/s00233-006-0667-5},
}
%%%%%
\bib{VV02}{article}{
   author={V\'azquez, Juan Luis},
   author={Vitillaro, Enzo},
   title={Heat equation with dynamical boundary conditions of reactive type},
   journal={Comm. Partial Differential Equations},
   volume={33},
   date={2008},
%   number={4-6},
   pages={561--612},
%   issn={0360-5302},
%   review={\MR{2424369}},
%   doi={10.1080/03605300801970960},
}
%%%%%
\bib{Vi}{article}{
   author={Vitillaro, Enzo},
   title={Global existence for the heat equation with nonlinear dynamical
   boundary conditions},
   journal={Proc. Roy. Soc. Edinburgh Sect. A},
   volume={135},
   date={2005},
%   number={1},
   pages={175--207},
%   issn={0308-2105},
%   review={\MR{2119848}},
%   doi={10.1017/S0308210500003838},
}
%%%%%
\bib{W}{article}{
   author={Weissler, Fred B.},
   title={Existence and nonexistence of global solutions for a semilinear
   heat equation},
   journal={Israel J. Math.},
   volume={38},
   date={1981},
%   number={1-2},
   pages={29--40},
%   issn={0021-2172},
%   review={\MR{0599472}},
%   doi={10.1007/BF02761845},
}
\end{biblist}
\end{bibdiv}
\end{document}